\documentclass[12pt]{amsart}
\usepackage{amsmath,amssymb,amscd,latexsym}

\newtheorem{theorem}{Theorem}[section]
\newtheorem{proposition}[theorem]{Proposition}

\newtheorem{lemma}[theorem]{Lemma}
\newtheorem{corollary}[theorem]{Corollary}
\newtheorem{remark}[theorem]{Remark}
\newtheorem{${}$}[theorem]{${}$}

\begin{document}

\title{Backward Uniqueness for Parabolic Operators\\ with Variable Coefficients in a Half Space}

\author{Jie Wu,\,\,\,\,\,\, Liqun Zhang$^{\ast}$}

\begin{figure}[b]

{\small
\begin{tabular}{ll}
$^{\ast}$The research is partially supported by the Chinese NSF under grant 10325104,
 the innovation\\ program at CAS and National Basic Research Program
of China under grant 2011CB808002.
\end{tabular}}
\end{figure}

\date{June, 2013}

\begin{abstract}It is shown that a function $u$ satisfying
$|\partial_tu+\sum_{i,j}\partial_i(a^{ij}\partial_ju)|\leq
N(|u|+|\nabla u|)$, $|u(x,t)|\leq Ne^{N|x|^2}$ in
$\mathbb{R}^n_+\times[0,T]$ and $u(x,0)=0$ in $\mathbb{R}^n_+$ under
certain conditions on $\{a^{ij}\}$ must vanish identically in
$\mathbb{R}^n_+\times[0,T]$. The main point of the result is that
the conditions imposed on $\{a^{ij}\}$ are of this type: $\{a^{ij}\}$
are Lipschitz and $|\nabla_xa^{ij}(x,t)|\leq \frac{E}{|x|}$, where
$E$ is less than a given number, and the conditions are optimal in some sense.
\end{abstract}

\maketitle

\\


\section{Introduction}

Let $U$ be a domain in $\mathbb{R}^n$ and $\emph{P}$ be a backward
parabolic operator on $U\times[0,T]$,
$$
P=\partial_t+\sum_{i,j}\partial_i(a^{ij}\partial_j)=\partial_t+\nabla\cdot(A\nabla),
$$
where $A(x,t)=(a^{ij}(x,t))^n_{i,j=1}$ is a real symmetric matrix such that for some $\Lambda\geq\lambda>0$,\\
\begin{equation}\label{a1}
\lambda|\xi|^2\leq \sum_{i,j}a^{ij}(x,t)\xi_i\xi_j\leq
\Lambda|\xi|^2, ~\forall \xi\in\mathbb{R}^n.
\end{equation}
Consider a function $u$ which satisfies
$$
\left\{
  \begin{array}{ll}
    |Pu|\leq N(|u|+|\nabla u|)~~~~&in~~ U\times[0,T], \\
    |u(x,t)|\leq Ne^{N|x|^2}~~~~&in~~ U\times[0,T], \\
    u(x,0)=0~~~~~~&in~~U,
  \end{array}
\right.
$$
where $N$ is a given positive constant.

The backward uniqueness (BU) problem is: does $u$ vanish identically in $U\times[0,T]$?
If so, we say that $U$ is a BU domain for the operator $P$.\\

We should point out that there is no boundary condition about $u$ on the boundary of the domain $U$. The BU problem is recently discussed by
L. Escauriaza, G. Seregin and V. \v{S}ver\'{a}k in \cite{ESS}. It appeared in many problems, for example, in the control theory for PDEs and the regularity theory of parabolic equations. Especially, it plays an important role in the regularity theory of the Navier-Stokes equations, see \cite{ESS3}.

When $P$ is the backward heat operator, there are many
results already on various domains, such as, on the exterior of a ball
$\mathbb{R}^n\backslash B_R$ \cite{ESS} and a half space
$\mathbb{R}^n_+=\{x\in\mathbb{R}^n| x_n>0\}$ \cite{ESS2,ESS3} by L. Escauriaza, G. Seregin and V. \v{S}ver\'{a}k, and on some cones by Lu Li and V. \v{S}ver\'{a}k \cite{LlS}. Those are all proved to be BU domains for
the backward heat operator. On the other hand, any bounded domain is not
BU domain, see \cite{Jone,Lit}.

When $P$ is in general, variable coefficients, there are few results
have been proved, while some related results have already been
obtained. In particular, L. Escauriaza and F. J. Fern\'{a}ndez
proved a unique continuation property when $\{{a^{ij}}\}$ are
Lipschitz in \cite{ESS4}. Then it implies immediately that if
$U\subset V$ and $U$ is a BU domain, so is $V$. Recently, Tu A.
Nguyen in \cite{Tu} proved a conjecture of E. M. Landis and O. A.
Oleinik \cite{LO} which implies that $\mathbb{R}^n$ and $\mathbb{R}^n_+$ are
BU domains under the conditions that $|\nabla_xa^{ij}(x,t)|$ and
$|\partial_ta^{ij}(x,t)|$ are bounded and the decay at infinity
conditions that
\begin{equation}\label{Tu}
|\nabla_xa^{ij}(x,t)|\leq M\langle x\rangle^{-1-\varepsilon}, ~~~~|a^{ij}(x,t)-a^{ij}(x,s)|
\leq M\langle x\rangle^{-1}|t-s|^{1/2},
\end{equation}
where $\langle x\rangle=\sqrt{1+|x|^2}$ and $\varepsilon>0$.

This paper can be regarded as a continuation of the above results.
Since $\mathbb{R}^n\backslash B(R)\subset\mathbb{R}^n$ and
$\mathbb{R}^n_+$ can be treated as a subset of
$\mathbb{R}^n\backslash B(R)$, by the unique continuation property,
we could only consider the case of $\mathbb{R}^n_+$. Also general
simply connected domains may be mapped onto $\mathbb{R}^n_+$. Then
we focus on operator $P$ with variable coefficients on the domain
$\mathbb{R}^n_+$. Our main result is the following.

\begin{theorem}\label{thm1}
Suppose $\{a^{ij}\}$ satisfy (\ref{a1}), and for some constants $E,M,N>0$,
\begin{equation}\label{a2}
|\nabla_xa^{ij}(x,t)|+|\partial_ta^{ij}(x,t)|\leq M, ~~\forall(x,t)\in \mathbb{R}^n_+\times[0,T],
\end{equation}
and
\begin{equation}\label{a3}
|\nabla_xa^{ij}(x,t)|\leq \frac{E}{|x|},  ~~\forall(x,t)\in \mathbb{R}^n_+\times[0,T].
\end{equation}
Assume that $u$ satisfies
\begin{eqnarray}\label{a4}
\left\{
  \begin{array}{ll}
    |Pu|\leq N(|u|+|\nabla u|)~~~~&in~~ \mathbb{R}^n_+\times[0,T],  \\
    |u(x,t)|\leq Ne^{N|x|^2}~~~~&in~~ \mathbb{R}^n_+\times[0,T],  \\
    u(x,0)=0~~~~&in~~\mathbb{R}^n_+.
  \end{array}
\right.
\end{eqnarray}
Then there exists a constant $E_0=E_0(n,\Lambda,\lambda)$, such that
when $E<E_0$, $u(x,t)\equiv 0$ in $\mathbb{R}^n_+\times[0,T]$.
\end{theorem}

We remark that our assumptions are optimal in some sense. From the
counterexamples constructed by A. Plis \cite{Pl},
K. Miller \cite{Mi} and N. Mandache \cite{Man}, we can see that to
ensure BU, certain regularity of the coefficients should be
required. Moreover, Tu A. Nguyen proved in \cite{Tu} that the
regularity conditions (\ref{a2}) and the decay at infinity
conditions (\ref{Tu}) will be sufficient to ensure BU. However, are
the decay at infinity conditions necessary? Or, are conditions
(\ref{a2}) alone enough to guarantee BU? Here, we show that
conditions (\ref{a2}) are not enough and the decay at infinity
conditions (\ref{a3}) in Theorem \ref{thm1},
$|\nabla_xa^{ij}(x,t)|\leq \frac{E}{|x|}$, where $E$ is small, are optimal in some sense.

First, we copy the examples given by N. Mandache in \cite{Man}.

\begin{proposition}\label{Manexample}
There exist smooth functions $u$, $b_{11}$, $b_{12}$, $b_{22}$ and
continuous functions $d_1$, $d_2$ defined on $\mathbb{R}^3\ni(s,x,y)$, with the following properties:\\
i) $u$ is the solution of the equation
\begin{equation}\label{Man1}
\partial_s^2u+\partial_x((b_{11}+d_1)\partial_xu)+\partial_y(b_{12}\partial_xu)
+\partial_x(b_{12}\partial_yu)+\partial_y((b_{22}+d_2)\partial_yu)=0.
\end{equation}
ii) There is a $T>0$ such that $supp~u=(-\infty,T]\times\mathbb{R}^2$.\\
iii) $u$, $b_{ij}$ and $d_i$ are periodic in x and in y with period $2\pi$.\\
iv) $d_1$ and $d_2$ do not depend on $x$ and $y$ and are H\"{o}lder continuous of order $\alpha$ for
all $\alpha<1$.\\
v) $\frac{1}{2}<\left(
                                      \begin{array}{cc}
                                        d_1+b_{11} & b_{12} \\
                                        b_{12} & d_2+b_{22} \\
                                      \end{array}
                                    \right)<2
$ on $\mathbb{R}^3$.\\

Furthermore, there are also functions as above, satisfying conditions i)-v) except
that (\ref{Man1}) is replaced with the parabolic equation:

\begin{equation}\label{Man2}
\partial_su=\partial_x((b_{11}+d_1)\partial_xu)+\partial_y(b_{12}\partial_xu)
+\partial_x(b_{12}\partial_yu)+\partial_y((b_{22}+d_2)\partial_yu).
\end{equation}
\end{proposition}

The solution of (\ref{Man2}) implies that the H\"{o}lder regularity
in the time variable is not enough for BU.
Hence it is reasonable for us to assume that $|\partial_ta^{ij}(x,t)|$ are bounded in Theorem \ref{thm1}.\\

Next we consider the requirement of the regularity in the space
variable. Assume that $u$ is the solution of (\ref{Man1}). We denote
that
$$v(t,s,x,y)=u(T+s+t,x,y);$$
$$\bar{b_{ij}}(t,s,x,y)=b_{ij}(T+s+t,x,y);$$
$$\bar{d_i}(t,s)=d_i(T+s+t).$$
in
$$[-1,0]\times\mathbb{R}^3_+=\{(t,s,x,y)|t\in[-1,0],s>0,x\in\mathbb{R},y\in\mathbb{R}\}.$$
Then
$$
\partial_tv-[\partial_s^2v+\partial_x((\bar{b_{11}}+\bar{d_1})\partial_xv)+\partial_y(\bar{b_{12}}\partial_xv)
+\partial_x(\bar{b_{12}}\partial_yv)+\partial_y((\bar{b_{22}}+\bar{d_2})\partial_yv)]
-\partial_sv=0.
$$
By $ii)$ of Proposition \ref{Manexample}, $v(0,s,x,y)=0$ and $v$ is
nonzero in $[-1,0]\times\mathbb{R}^3_+$, thus BU fails. It shows
that the H\"{o}lder regularity in the space variable is not enough
for BU,
hence it is also reasonable for us to assume that $|\nabla_xa^{ij}(x,t)|$
are bounded in Theorem \ref{thm1}.\\

Now we consider the decay at infinity conditions. We could construct
an example as follows. Consider a cone $\mathcal{C}_{\theta_0}$ with
opening angle $\theta_0$ and the system
\begin{eqnarray}\label{CE1}
\left\{
  \begin{array}{ll}
    \partial_tu+\Delta u=0~~~~&in~~ \mathcal{C}_{\theta_0}\times[0,T],  \\
    |u(x,t)|\leq N~~~~&in~~ \mathcal{C}_{\theta_0}\times[0,T], \\
    u(x,0)=0~~~~~&in~~\mathcal{C}_{\theta_0}.
  \end{array}
\right.
\end{eqnarray}

In \cite{LlS}, L. Escauriaza gave an example to show that the above system has a nonzero solution
when $\theta_0<\frac{\pi}{2}$ and the result of Lu Li and V. \v{S}ver\'{a}k implied that the system has only zero solution
when $\pi>\theta_0>2\arccos(1/\sqrt{3})\approx109.5^\circ$. Now we consider  a cone of dimension 2,
$$\mathcal{C}_{\theta_0}=\{(r,\theta)| 0<\theta<\theta_0\}, ~~ (0<\theta_0<\pi)$$
and $u(x_1,x_2,t)$ is the solution of system (\ref{CE1}) in
dimension 2, where $$\left\{
      \begin{array}{ll}
        x_1=r\cos\theta & \hbox{} \\
        x_2=r\sin\theta. & \hbox{}
      \end{array}
    \right.$$
Let
$$\bar{\theta}=l\theta, ~~\text{with}~~l=\frac{\pi}{\theta_0}>1,$$
and
$$
\left\{
      \begin{array}{ll}
        y_1=r\cos\bar{\theta} & \hbox{} \\
        y_2=r\sin\bar{\theta}, & \hbox{}
      \end{array}
    \right.
$$
then $(y_1,y_2)\in \mathbb{R}^2_+$. We denote
$$v(y_1,y_2,t)=u(x_1,x_2,t),~~(y_1,y_2,t)\in\mathbb{R}^2_+\times[0,T].$$
By simple calculation,
\begin{eqnarray}\label{CE2}
\begin{split}
\partial_{x_1}^2+\partial_{x_2}^2=&\partial_{r}^2+\frac{\partial_r}{r}+\frac{\partial_{\theta}^2}{r^2}
=\partial_{r}^2+\frac{\partial_r}{r}+l^2\frac{\partial_{\bar{\theta}}^2}{r^2}\\
=&(\partial_{r}^2+\frac{\partial_r}{r}+\frac{\partial_{\bar{\theta}}^2}{r^2})+
(l^2-1)\frac{\partial_{\bar{\theta}}^2}{r^2}\\
=&\partial_{y_1}^2+\partial_{y_2}^2+(l^2-1)(\frac{y_2^2}{r^2}\partial_{y_1}^2+\frac{y_1^2}{r^2}\partial_{y_2}^2
-2\frac{y_1y_2}{r^2}\partial_{y_1 y_2}-\frac{y_1}{r^2}\partial_{y_1}-\frac{y_2}{r^2}\partial_{y_2})\\
=&[1+(l^2-1)\frac{y_2^2}{r^2}]\partial_{y_1}^2+[1+(l^2-1)\frac{y_1^2}{r^2}]\partial_{y_2}^2\\
&-2(l^2-1)\frac{y_1y_2}{r^2}\partial_{y_1 y_2}-(l^2-1)(\frac{y_1}{r^2}\partial_{y_1}+\frac{y_2}{r^2}\partial_{y_2})\\
=&\nabla\cdot(A\nabla v),
\end{split}
\end{eqnarray}
where
$$A(y_1,y_2)=\left(
                       \begin{array}{cc}
                         1+(l^2-1)\frac{y_2^2}{r^2} & -(l^2-1)\frac{y_1y_2}{r^2} \\
                         -(l^2-1)\frac{y_1y_2}{r^2} & 1+(l^2-1)\frac{y_1^2}{r^2} \\
                       \end{array}
                     \right).
$$
Together with the equation in (\ref{CE1}) we can deduce that
\begin{equation}\label{CE3}
\partial_tv+\nabla\cdot(A\nabla v)=0,
\end{equation}
and $A$ is positive since $l>1$.

Denote
$$w(y_1,y_2,t)=v(y_1,y_2+1,t), ~~(y_1,y_2,t)\in\mathbb{R}^2_+\times[0,T],$$
and
$$
B(y_1,y_2)=A(y_1,y_2+1)\equiv\left(
                               \begin{array}{cc}
                                 b^{11} & b^{12} \\
                                 b^{21} & b^{22} \\
                               \end{array}
                             \right).
$$
Direct calculations give us
$$
|\nabla b^{ij}|\leq l^2-1 ~~and ~~|\nabla b^{ij}|\leq \frac{l^2-1}{r}\equiv\frac{E_1}{r}.
$$
By (\ref{CE3}) we have
$$\partial_tw+\nabla\cdot(B\nabla w)=0.$$

By the notations of $u$ and $w$, we see that $w$ is a solution of
the following system:
\begin{eqnarray}\label{CE4}
\left\{
  \begin{array}{ll}
    \partial_tw+\nabla\cdot(B\nabla w)=0~~~~&in~~ \mathbb{R}^2_+\times[0,T],  \\
    |w(y,t)|\leq N~~~~~~~~&in~~ \mathbb{R}^2_+\times[0,T], \\
    w(y,0)=0~~~~~~~~&in~~\mathbb{R}^2_+.
  \end{array}
\right.
\end{eqnarray}

By the result of Li and \v{S}ver\'{a}k \cite{LlS}, we conclude that
when $$E_1<(\frac{\pi}{2\arccos(1/\sqrt{3})})^2-1,$$ we have
$$1<l<\frac{\pi}{2\arccos(1/\sqrt{3})},\quad
2\arccos(1/\sqrt{3})<\theta_0<\pi,$$ and then $u\equiv0$ and thus
$w\equiv0$.\\
When $E_1>3$, we have $l>2$, $\theta_0<\frac{\pi}{2}$,
and then (\ref{CE1}) has a nonzero solution and thus (\ref{CE4})
must also has a nonzero solution $w$. Otherwise, if $w\equiv0$, this means that $u=0$ in $D\times[0, T]$, where $D$ is a certain subregion of $\mathcal{C}_{\theta_0}$. Then by the unique continuation result, we have $u\equiv0$, which contradicts that $u$ is nonzero. In this case BU fails, although $|\nabla
b^{ij}|$ are bounded and $|\nabla b^{ij}|\leq \frac{E_1}{r}$.

The example above shows that the decay at infinity conditions, those are assumptions in (4), where $E$ is less than a given constant, are optimal in some sense.

To prove Theorem \ref{thm1} we need to obtain the
corresponding Carleman inequalities. Now we introduce two Carleman
inequalities for the case of variable coefficients. They are
generalizations of the two Carleman inequalities for the case of
constant coefficients, as shown in \cite{ESS2,ESS3}.

\begin{proposition}\label{Prop-C1}
Suppose $\{a^{ij}\}$ satisfy (\ref{a1}) and
$$|\nabla_xa^{ij}(x,t)|+|\partial_ta^{ij}(x,t)|\leq M,~~|\nabla_xa^{ij}(x,t)|\leq
\frac{E}{|x|},~~~\forall(x,t)\in \mathbb{R}^n\times(0,2).$$
Then there exists a constant $K=K(n,\Lambda,\lambda,M,E)$, such that
for any $u\in C^\infty_0(\mathbb{R}^n\times(0,2))$ and any number
$\gamma>0$,
\begin{equation}\label{C1}
\begin{split}
\int_{\mathbb{R}^n\times (0,2)}e^{2\gamma(t^{-K}-1)-\frac{b|x|^2+K}{t}}(|u|^2+|\nabla u|^2)dxdt\\
\leq \int_{\mathbb{R}^n\times (0,2)}e^{2\gamma(t^{-K}-1)-\frac{b|x|^2+K}{t}}|Pu|^2dxdt,
\end{split}
\end{equation}
where $b=\frac{1}{8\Lambda}$.
\end{proposition}

\begin{proposition}\label{Prop-C2}
Suppose $\{a^{ij}\}$ satisfy (\ref{a1}) and
$$|\nabla_xa^{ij}(x,t)|+|\partial_ta^{ij}(x,t)|\leq M,~~|\nabla_xa^{ij}(x,t)|\leq
\frac{E}{|x|},~~~\forall(x,t)\in \mathbb{R}^n_+\times(0,1).$$
Let $Q=\{(x,t)|x_n\geq 1, t\in(0,1)\}$ and
$$\psi(x)=|x|^2-2\frac{\Lambda}{\lambda}|x|x_n+2(\frac{\Lambda}{\lambda})^2x_n^2.$$
Then there exist positive constants $E_0=E_0(n, \Lambda, \lambda)$,
$\alpha=\alpha(n,\Lambda,\lambda,E)\in(1,2)$, $b=b(\Lambda,\lambda)$
and $K=K(n,\Lambda,\lambda,M,E)$ such that when $E< E_0$, for any
function $u\in C^\infty_0(Q)$ and any number $\gamma>0$, we have
\begin{equation}\label{C2}
\begin{split}
\int_{Q}e^{2\gamma (t^{-K}-1)x_n^{\alpha}-\frac{b\psi(x)+K}{t}}(|u|^2+|\nabla u|^2)dxdt\\
\leq \int_{Q}e^{2\gamma
(t^{-K}-1)x_n^{\alpha}-\frac{b\psi(x)+K}{t}}|Pu|^2dxdt.
\end{split}
\end{equation}

\end{proposition}
\begin{remark}
In fact, we can take
$E_0=\frac{\lambda}{16n^2\frac{\Lambda}{\lambda}(\frac{\Lambda}{\lambda}+1)}$,
$\alpha=1+\frac{E}{E_0}$ and
$b=\frac{1}{64\Lambda(\frac{\Lambda}{\lambda}+1)^4}$ in Proposition
\ref{Prop-C2} which can be seen from the proof.
\end{remark}

Carleman inequality (\ref{C2}) is the key results in this paper.
Assuming it, there is only a standard argument by following the
corresponding parts of Escauriaza, Seregin, and \v{S}ver\'{a}k in
\cite{ESS,ESS2} to prove Theorem \ref{thm1}. In the
establishment of Carleman inequality (\ref{C2}), the construction of
the function $\psi$ is crucial.

\begin{remark}
It is worthwhile to note that Carleman inequality (\ref{C1}) does
not require the smallness of $E$, while Carleman inequality
(\ref{C2}) does, which is stronger.
\end{remark}

Moreover, the Carleman inequality (\ref{C2}) for the parabolic
operators with variable coefficients in a half space is stronger
than the one for the case of constant coefficients, as shown in
\cite{ESS2,ESS3}. When $P$ is the backward heat operator, there are
 two Carleman inequalities to prove BU. The first one implies an
exponential decay of the solution, which enable us to apply
the second one to prove BU. And here, we just need one Carleman inequality (\ref{C2})
to prove BU. We list Carleman inequality (\ref{C1}) here just for
comparison with the case of constant coefficients.

The paper organized as follows. We first make use of Carleman
inequality (\ref{C2}) to prove Theorem \ref{thm1} in next section. Then we
prove the two Carleman inequalities Proposition \ref{Prop-C1} and
Proposition \ref{Prop-C2} in the last section.

\section{Proof of Theorem \ref{thm1}}

In this section, we prove the main theorem by assuming Proposition
1.4 first. Then we shall prove the Carleman inequalities in next
section.

We always assume that $T=1$ and extend $u$ and $a^{ij}$ by the following way:
$$u(x,t)=0,~if~t<0;$$
$$a^{ij}(x,t)=a^{ij}(x,0),~if~t<0.$$
We denote $e_n=(0,0,...,0,1)$.

The next lemma implies Theorem \ref{thm1} immediately.
\begin{lemma}
Suppose $\{a^{ij}\}$ and $u$ satisfy assumptions (\ref{a1}),
(\ref{a2})-(\ref{a4}). Then there exists
$T_1=T_1(\Lambda,\lambda,N)\in(0,\frac{1}{2})$, such that
$$u(x,t)\equiv0$$
in $\mathbb{R}^n_+\times(0,T_1)$.
\end{lemma}
\emph{Proof}. We make use of Carleman inequality (\ref{C2}) to prove
this lemma. We mainly follow the arguments of corresponding parts of
Escauriaza, Seregin and \v{S}ver\'{a}k in \cite{ESS,ESS2}. By the regularity theory for solutions
of parabolic equations, we have
\begin{equation}\label{pt1}
|u(x,t)|+|\nabla u(x,t)|\leq C(n,\Lambda,\lambda,M,N)e^{2N|x|^2}
\end{equation}
for $(x,t)\in (\mathbb{R}^n_++e_n)\times(0,\frac{1}{2})$.
Let
\begin{equation}\label{T_1variable}
T_1=\min\{\frac{b}{32N},\frac{1}{12N^2},\frac{1}{2}\},
\end{equation}
where $b$ is the one in Proposition \ref{Prop-C2}. Let$\tau=\sqrt{2T_1}$. \\
We denote
$$v(y,s)=u(\tau y,\tau^2s-T_1)$$
and
$$\tilde{a}^{ij}(y,s)={a}^{ij}(\tau y,\tau^2s-T_1)$$
for $(y,s)\in\mathbb{R}^n_+\times(0,1)$. Then it is easy to see
$$|\nabla_y\tilde{a}^{ij}(y,s)|+|\partial_s\tilde{a}^{ij}(y,s)|\leq \tau M\leq M,$$
and
$$|\nabla_y\tilde{a}^{ij}(y,s)|=\tau|\nabla a^{ij}(\tau y, \tau^2s-T_1)|\leq\tau\frac{E}{|\tau y|}=\frac{E}{|y|}.$$
We denote
$$\tilde{P}v=\partial_sv+\sum_{ij}\partial_{i,y}(\tilde{a}^{ij}\partial_{j,y}v),$$
by our notation and (\ref{a4}),
\begin{equation}\label{pt2}
|\tilde{P}v|\leq\tau N(|v|+|\nabla v|),
\end{equation}
for $(y,s)\in\mathbb{R}^n_+\times(0,1)$. From (\ref{pt1}), we have
\begin{equation}\label{pt3}
|v(y,s)|+|\nabla v(y,s)|\leq C(n,\Lambda,\lambda,M,N)e^{2N\tau^2|y|^2};
\end{equation}
for $(y,s)\in (\mathbb{R}^n_++\frac{1}{\tau}e_n)\times(0,1)$; and
\begin{equation}\label{pt4}
v(y,s)=0,
\end{equation}
for $(y,s)\in \mathbb{R}^n_+\times(0,\frac{1}{2}]$.

In order to apply Carleman inequality (\ref{C2}), we choose two
smooth cut-off functions such that
$$
\eta_1(p)=\left\{
                 \begin{array}{ll}
                   0, & \hbox{if $p<\frac{1}{\tau}+1$;} \\
                   1, & \hbox{if $p>\frac{1}{\tau}+2$.}
                 \end{array}
               \right.
$$
And
$$
\eta_2(q)=\left\{
            \begin{array}{ll}
              0, & \hbox{if $q<-\frac{3}{4}$;} \\
              1, & \hbox{if $q>-\frac{1}{2}$.}
            \end{array}
          \right.
$$
All functions take values in $[0,1]$ and $|\eta_1'|$, $|\eta_1''|$,
$|\eta'_2|$ and $|\eta''_2|$ are all bounded. Denote
$$f(s)=s^{-K}-1$$ and
$$C_\star=1+\sup_{\frac{1}{2}<s<1 \atop \frac{1}{\tau}+1<y_n<\frac{1}{\tau}+2}\{f(s)y_n^\alpha\}
=1+f(\frac{1}{2})(\frac{1}{\tau}+2)^\alpha.$$ Set
$$\eta(y,s)=\eta_1(y_n)\eta_2(\frac{f(s)y_n^\alpha}{2C_\star}-1),$$
and $w=\eta v$. Then $supp~w\subset Q$, and
\begin{eqnarray*}
\begin{split}
|\tilde{P}w|=&|\eta \tilde{P}v+v\tilde{P}\eta+2\tilde{a}^{ij}\partial_i\eta\partial_j v|\\
\leq&|\eta \tilde{P}v|+C(n,\Lambda,M)\chi(|v|+|\nabla v|)(|\partial_s\eta|+|\nabla\eta|+|\nabla^2\eta|),
\end{split}
\end{eqnarray*}
where $\chi$ is the characteristic function of the set
$$\Omega=\{(y,s)|\frac{1}{2}<s<1,0<\eta<1\}.$$
By (\ref{pt2}), we have
\begin{eqnarray*}
\begin{split}
|\tilde{P}w|\leq&\eta\tau N(|v|+|\nabla v|)+C(n,\Lambda,M)\chi(|v|+|\nabla v|)(|\partial_s\eta|+
|\nabla\eta|+|\nabla^2\eta|)\\
\leq&\tau N(|w|+|\nabla w|)+C(n,\Lambda,M,N)\chi(|v|+|\nabla v|)(|\partial_s\eta|+|\nabla\eta|+|\nabla^2\eta|).
\end{split}
\end{eqnarray*}
Notice that $\frac{1}{2}<s<1$ in $\Omega$, and when $\frac{1}{2}<s<1$,
$$|\partial_s\eta|+|\nabla\eta|+|\nabla^2\eta|\leq C(n,\Lambda,\lambda,M,E)y_n^{\alpha}\leq
C(n,\Lambda,\lambda,M,E)y_n^{2},$$
then
\begin{equation}\label{pt5}
|\tilde{P}w|\leq \tau N(|w|+|\nabla w|)+C(n,\Lambda,\lambda,M,E,N)\chi(|v|+|\nabla v|)y_n^2.
\end{equation}
Moreover,
\begin{eqnarray*}
\begin{split}
\Omega=&\{(y,s)|\frac{1}{2}<s<1,\eta_1>0,0<\eta_2<1\}\\
&\bigcup\{(y,s)|\frac{1}{2}<s<1,0<\eta_1<1,\eta_2=1\}\\
=&\{(y,s)|\frac{1}{2}<s<1,y_n>\frac{1}{\tau}+1,\frac{1}{2}<\frac{f(s)y_n^\alpha}{C_\star}<1\}\\
&\bigcup\{(y,s)|\frac{1}{2}<s<1,\frac{1}{\tau}+1<y_n<\frac{1}{\tau}+2,\frac{f(s)y_n^\alpha}{C_\star}\geq1\}.
\end{split}
\end{eqnarray*}
By the choice of $C_\star$ we obtain that the second set of the right side of the above identity is empty, then
\begin{equation}\label{pt6}
\Omega=\{(y,s)|\frac{1}{2}<s<1,y_n>\frac{1}{\tau}+1,\frac{1}{2}<\frac{f(s)y_n^\alpha}{C_\star}<1\}.
\end{equation}
By (\ref{pt3}), in the support of $w$ we have
$$e^{2\gamma f(s)y_n^\alpha-\frac{b\psi(y)+K}{s}}(|v|+|\nabla v|)^2\leq Ce^{2\gamma f(s)y_n^\alpha-\frac{b\psi(y)+K}{s}+4N\tau^2|y|^2}.$$
Notice that $\psi(y)\geq\frac{|y|^2}{2}$ and $\tau=\sqrt{2T_1}$, then we have
$$e^{2\gamma f(s)y_n^\alpha-\frac{b\psi(y)+K}{s}}(|v|+|\nabla v|)^2\leq
Ce^{2\gamma f(s)y_n^\alpha-\frac{b}{2}|y|^2+8NT_1|y|^2}.$$ By
(\ref{T_1variable}), we know that $T_1\leq\frac{b}{32N}$, then in
$supp~w$ we have
\begin{equation}\label{pt7}
e^{2\gamma f(s)y_n^\alpha-\frac{b\psi(y)+K}{s}}(|v|+|\nabla v|)^2
\leq Ce^{2\gamma f(s)y_n^\alpha-\frac{b}{4}|y|^2}.
\end{equation}
Although $supp~w$ may be unbounded, $supp~w\subset Q$
and (\ref{pt7}) allow us to claim the validity of Proposition
\ref{Prop-C2} for $w$. Then by  Carleman inequality (\ref{C2}),
together with (\ref{pt5}), we have
\begin{eqnarray*}
\begin{split}
J&\equiv\int_Qe^{2\gamma f(s)y_n^\alpha-\frac{b\psi(y)+K}{s}}(|w|^2+|\nabla w|^2)dyds\\
&\leq \int_Qe^{2\gamma f(s)y_n^\alpha-\frac{b\psi(y)+K}{s}}|\tilde{P}w|^2dyds\\
&\leq 3\tau^2N^2J+C\int_Qe^{2\gamma f(s)y_n^\alpha-\frac{b\psi(y)+K}{s}}\chi(|v|+|\nabla v|)^2y_n^4dyds.
\end{split}
\end{eqnarray*}
By (\ref{T_1variable}), we know that $3\tau^2N^2=6T_1
N^2\leq\frac{1}{2}$. We deduce from the above inequality that
$$J\leq C\int_\Omega e^{2\gamma f(s)y_n^\alpha-\frac{b\psi(y)+K}{s}}(|v|+|\nabla v|)^2y_n^4dyds.$$
By (\ref{pt7}), we have
$$J\leq C\int_\Omega e^{2\gamma f(s)y_n^\alpha-\frac{b}{4}|y|^2}y_n^4dyds.$$
By (\ref{pt6}) we have that $f(s)y_n^\alpha<C_\star$ in $\Omega$, then
\begin{equation}\label{pt8}
J\leq C e^{2\gamma C_\star}\int_\Omega e^{-\frac{b}{4}|y|^2}y_n^4dyds\leq C e^{2\gamma C_\star}.
\end{equation}
On the other hand, we denote
$$\Omega_1=\{(y,s)|0<s<1,\eta=1\}=\{(y,s)|0<s<1, y_n\geq\frac{1}{\tau}+2, \frac{f(s)y_n^\alpha}{C_\star}\geq1\},$$
\begin{equation}\label{omega2variable}
\Omega_2=\{(y,s)|0<s<1, y_n\geq\frac{1}{\tau}+2, \frac{f(s)y_n^\alpha}{C_\star}\geq2\}.
\end{equation}
Obviously $\Omega_2\subset\Omega_1$ and $w=v$ in $\Omega_1$. Then
\begin{eqnarray*}
\begin{split}
J\geq&\int_{\Omega_1}e^{2\gamma f(s)y_n^\alpha-\frac{b\psi(y)+K}{s}}(|v|^2+|\nabla v|^2)dyds\\
\geq&\int_{\Omega_2}e^{2\gamma f(s)y_n^\alpha-\frac{b\psi(y)+K}{s}}(|v|^2+|\nabla v|^2)dyds\\
\end{split}
\end{eqnarray*}
By (\ref{omega2variable}), we know that $f(s)y_n^\alpha\geq2C_\star$
in $\Omega_2$. Hence
\begin{equation}\label{pt9}
J\geq e^{4\gamma C_\star}\int_{\Omega_2}e^{-\frac{b\psi(y)+K}{s}}(|v|^2+|\nabla v|^2)dyds.
\end{equation}
Combining (\ref{pt8}) and (\ref{pt9}), we have
$$\int_{\Omega_2}e^{-\frac{b\psi(y)+K}{s}}(|v|^2+|\nabla v|^2)dyds\leq Ce^{-2\gamma C_\star}.$$
Passing to the limit as $\gamma\rightarrow+\infty$, we obtain
$v(y,s)=0$ in $\Omega_2$. Using unique continuation though spatial
boundaries (see \cite{ESS4}), we obtain that $v(y,s)\equiv 0$ in
$\mathbb{R}^n_+\times(0,1)$. That is, $u(x,t)\equiv 0$ in
$\mathbb{R}^n_+\times(0,T_1)$. Thus we proved this lemma.

\section{Proof of Carleman Inequalities}

In this section, we shall prove two Carleman Inequalities which is
the crucial part of the whole argument. The main idea is to choose a
proper weight function $G$. We denote
$$\tilde{\Delta}u=\partial_i(a^{ij}\partial_ju).$$ Here and in the
following argument, we use the summation convention on the repeated
indices. We shall make use of the following lemma which is due to L.
Escauriaza and F. J. Fern\'{a}ndez in \cite{ESS4} (see also
\cite{Tu}).
\begin{lemma}
Suppose $\sigma(t):\mathbb{R}_+\rightarrow \mathbb{R}_+$ is a smooth
function, $\alpha$ is a real number, $F$ and $G$ are differentiable
functions and $G>0$. Then the following identity holds for any $u\in
C^\infty_0(\mathbb{R}^n\times(0,T))$
\begin{equation}\label{pc0.1}
\begin{split}
&2\int_{\mathbb{R}^n\times(0,T)}\frac{\sigma^{1-\alpha}}{\sigma'}|Lu|^2Gdxdt
+\frac{1}{2}\int_{\mathbb{R}^n\times(0,T)}\frac{\sigma^{1-\alpha}}{\sigma'}u^2MGdsdt\\
&+\int_{\mathbb{R}^n\times(0,T)}\frac{\sigma^{1-\alpha}}{\sigma'}\langle A\nabla u,\nabla u\rangle
[(\log{\frac{\sigma}{\sigma'}})'+\frac{\partial_tG-\tilde{\Delta}G}{G}-F]Gdxdt\\
&+2\int_{\mathbb{R}^n\times(0,T)}\frac{\sigma^{1-\alpha}}{\sigma'}\langle D_G\nabla u,\nabla u\rangle Gdxdt
-\int_{\mathbb{R}^n\times(0,T)}\frac{\sigma^{1-\alpha}}{\sigma'}u\langle A\nabla u,\nabla F\rangle Gdxdt\\
&=2\int_{\mathbb{R}^n\times(0,T)}\frac{\sigma^{1-\alpha}}{\sigma'}LuPuGdxdt
\end{split}
\end{equation}
where
$$Lu=\partial_tu-\langle A\nabla{logG},\nabla u\rangle+\frac{Fu}{2}-\frac{\alpha\sigma'}{2\sigma}u,$$
$$M=(log{\frac{\sigma}{\sigma'}})'F+\partial_tF+(F-\alpha\frac{\sigma'}{\sigma})(\frac{\partial_tG-
\tilde{\Delta}G}{G}-F)-\langle A\nabla F,\nabla{logG}\rangle,$$
and
$$D^{ij}_G=a^{ik}\partial_{kl}(logG)a^{lj}
+\frac{\partial_l(logG)}{2}(a^{ki}\partial_ka^{lj}+a^{kj}\partial_ka^{li}-a^{kl}\partial_ka^{ij})+
\frac{1}{2}\partial_ta^{ij}.
$$
\end{lemma}
We first give a modification of this lemma which will be used in our
proof. Letting $\alpha=0$ and $\sigma(t)=e^t$ in Lemma 3.1, we
obtain the following identity for $u\in
C^\infty_0(\mathbb{R}^n\times(0,T))$
\begin{eqnarray*}
\begin{split}
&\frac{1}{2}\int_{\mathbb{R}^n\times(0,T)}u^2MGdxdt+\int_{\mathbb{R}^n\times(0,T)}[2\langle D_G\nabla u,\nabla u\rangle+\langle A\nabla u,\nabla u\rangle(\frac{\partial_tG-\tilde{\Delta}G}{G}-F)]Gdxdt\\
&-\int_{\mathbb{R}^n\times(0,T)}u\langle A\nabla u,\nabla F\rangle Gdxdt=2\int_{\mathbb{R}^n\times(0,T)}Lu(Pu-Lu)Gdxdt.
\end{split}
\end{eqnarray*}
If $\nabla F$ is differentiable, we can integrate by parts to obtain
\begin{eqnarray*}
\begin{split}
&-\int_{\mathbb{R}^n\times(0,T)}u\langle A\nabla u,\nabla F\rangle Gdxdt\\
=&\frac{1}{2}\int_{\mathbb{R}^n\times(0,T)}u^2\tilde{\Delta} FGdxdt+\frac{1}{2}\int_{\mathbb{R}^n\times(0,T)}u^2\langle A\nabla F,\nabla logG\rangle Gdxdt.
\end{split}
\end{eqnarray*}
The function $\nabla F$ may not be differentiable, so we approximate $F$ by some smooth function $F_0$ and use the above identity with $F_0$ in place of $F$, following Tu's idea in \cite{Tu}.  Thus we obtain the following result.

\begin{corollary}\label{CIforgeneral}
Suppose $F$ and $G$ are differentiable functions and $G>0$. Then the
following identity holds for any $u\in
C^\infty_0(\mathbb{R}^n\times(0,T))$
\begin{equation}\label{pc0.2}
\begin{split}
&\frac{1}{2}\int_{\mathbb{R}^n\times(0,T)}u^2M_0Gdxdt+\int_{\mathbb{R}^n\times(0,T)}
[2\langle D_G\nabla u,\nabla u\rangle+\langle A\nabla u,\nabla u\rangle
(\frac{\partial_tG-\tilde{\Delta}G}{G}-F)]Gdxdt\\
&-\int_{\mathbb{R}^n\times(0,T)}u\langle A\nabla u,\nabla (F-F_0)\rangle Gdxdt=
2\int_{\mathbb{R}^n\times(0,T)}Lu(Pu-Lu)Gdxdt,
\end{split}
\end{equation}
where
$$Lu=\partial_tu-\langle A\nabla u,\nabla{logG}\rangle+\frac{Fu}{2},$$
$$M_0=\partial_tF+F(\frac{\partial_tG-\tilde{\Delta}G}{G}-F)+\tilde{\Delta}F_0-
\langle A\nabla (F-F_0),\nabla{logG}\rangle,$$
and
$$D^{ij}_G=a^{ik}\partial_{kl}(logG)a^{lj}
+\frac{\partial_l(logG)}{2}(a^{ki}\partial_ka^{lj}+a^{kj}\partial_ka^{li}-
a^{kl}\partial_ka^{ij})+\frac{1}{2}\partial_ta^{ij}.
$$
\end{corollary}

\subsection{Proof of Proposition \ref{Prop-C1}.}
We use identity (\ref{pc0.2}) to prove Proposition \ref{Prop-C1}. In (\ref{pc0.2}),
we let $$G=e^{2\gamma(t^{-K}-1)-\frac{b|x|^2+K}{t}},$$
then
$$\frac{\partial_tG-\tilde{\Delta}G}{G}=\frac{b|x|^2-4b^2a^{ij}x_ix_j+K}{t^2}+\frac{2ba^{ii}+
2b\partial_ka^{kl}x_l}{t} -2\gamma Kt^{-K-1}.$$ Let
$$F=\frac{b|x|^2-4b^2a^{ij}x_ix_j+K}{t^2}+\frac{2ba^{ii}}{t}-2\gamma
Kt^{-K-1}-d(\frac{1}{t}+1),$$ where $d$ is a positive constant to be
determined. Set
$$F_0=\frac{b|x|^2-4b^2a^{ij}_\epsilon x_ix_j+K}{t^2}+\frac{2ba^{ii}_\epsilon}{t}-
2\gamma Kt^{-K-1}-d(\frac{1}{t}+1),$$
where $$a^{ij}_\epsilon(x,t)=\int_{\mathbb{R}^n}a^{ij}(x-y,t)\phi_\epsilon(y)dy,$$
$\phi$ is a mollifier, and $\epsilon=\frac{1}{2}$. \\

We denote by $I_n$ the identity matrix of $\mathbb{R}^n$, $C$ are
generic constants depending on $n,\Lambda,\lambda,M$ and $E$ in the
following arguments. We need some estimates which we list in the
following lemma.
\begin{lemma}\label{estimates1}
Set $b=\frac{1}{8\Lambda}$ and $K=12d$. For
$d=d(n,\Lambda,\lambda,M,E)$ large enough, we have
\begin{equation}\label{E1}
2D_G+A(\frac{\partial_t G-\tilde{\Delta}G}{G}-F)\geq(\frac{1}{t}+1)I_n;
\end{equation}
\begin{equation}\label{E2}
\partial_tF+F(\frac{\partial_tG-\tilde{\Delta}G}{G}-F)\geq \frac{db(|x|^2+1)}{4t^3};
\end{equation}
\begin{equation}\label{E3}
\tilde{\triangle}F_0\geq-\frac{C(|x|^2+1)}{t^3};
\end{equation}
\begin{equation}\label{E4}
|\nabla(F-F_0)|\leq\frac{C(|x|+1)}{t^2}.
\end{equation}
\end{lemma}
We will prove this lemma later.

First by applying (\ref{E1}) we have
\begin{equation}\label{pc1.1}
\begin{split}
&\int_{\mathbb{R}^n\times(0,2)}[2\langle D_G\nabla u,\nabla u\rangle+
\langle A\nabla u,\nabla u\rangle(\frac{\partial_tG-\tilde{\Delta}G}{G}-F)]Gdxdt\\
\geq&\int_{\mathbb{R}^n\times(0,2)}(\frac{1}{t}+1)|\nabla u|^2Gdxdt.
\end{split}
\end{equation}
Next we estimate $M_0$. By applying (\ref{E4}) we have
\begin{equation}\label{pc1.5}
|\langle A\nabla(F-F_0),\nabla logG\rangle|\leq\Lambda|\nabla(F-F_0)||\nabla logG|\leq
\frac{C(|x|+1)|x|}{t^3}\leq\frac{C(|x|^2+1)}{t^3}.
\end{equation}
Then by (\ref{E2}), (\ref{E3}) and (\ref{pc1.5}) we have
\begin{eqnarray*}
\begin{split}
M_0&=\partial_tF+F(\frac{\partial_tG-\tilde{\Delta}G}{G}-F)+\tilde{\Delta}F_0-
\langle A\nabla (F-F_0),\nabla{logG}\rangle\\
&\geq (\frac{db}{4}-C)\frac{|x|^2+1}{t^3},
\end{split}
\end{eqnarray*}
thus
\begin{equation}\label{pc1.6}
\frac{1}{2}\int_{\mathbb{R}^n\times(0,2)}u^2M_0Gdxdt
\geq (\frac{db}{8}-C)\int_{\mathbb{R}^n\times(0,2)}\frac{|x|^2+1}{t^3}u^2Gdxdt.
\end{equation}
By (\ref{E4}) and the Cauchy inequality we have
\begin{equation}\label{pc1.7}
\begin{split}
&|\int_{\mathbb{R}^n\times(0,2)}u\langle A\nabla u,\nabla (F-F_0)\rangle Gdxdt|\\
\leq& \Lambda\int_{\mathbb{R}^n\times(0,2)}|\nabla(F-F_0)||u||\nabla u|Gdxdt\\
\leq& C\int_{\mathbb{R}^n\times(0,2)}\frac{|x|+1}{t^2}|u||\nabla u|Gdxdt\\
\leq& C\int_{\mathbb{R}^n\times(0,2)}\frac{|x|^2+1}{t^3}|u|^2Gdxdt
+\int_{\mathbb{R}^n\times(0,2)}\frac{|\nabla u|^2}{t}Gdxdt.
\end{split}
\end{equation}
Finally, by (\ref{pc0.2}), (\ref{pc1.1}), (\ref{pc1.6}), (\ref{pc1.7}) and the Cauchy inequality, we have
$$
\int_{\mathbb{R}^n\times(0,2)}|Pu|^2Gdxdt
\geq (\frac{db}{8}-C)\int_{\mathbb{R}^n\times(0,2)}\frac{|x|^2+1}{t^3}|u|^2Gdxdt+\int_{\mathbb{R}^n\times(0,2)}|\nabla u|^2Gdxdt,\\
$$
if we choose $d$ large enough, we obtain
$$
\int_{\mathbb{R}^n\times(0,2)}|Pu|^2Gdxdt
\geq \int_{\mathbb{R}^n\times(0,2)}(|u|^2+|\nabla u|^2)Gdxdt.
$$
Thus we proved Carleman inequality (\ref{C1}).\\

\emph{Proof of Lemma \ref{estimates1}.} We estimate them one by one.\\

\emph{Estimate of $2D_G+A(\frac{\partial_t G-\tilde{\Delta}G}{G}-F)$.}\\

By direct calculations we have
\begin{eqnarray*}
\begin{split}
&2D_G+A(\frac{\partial_t G-\tilde{\Delta}G}{G}-F)\\
=&-\frac{4b}{t}A^2-\frac{2bx_l}{t}(a^{ki}\partial_ka^{lj}+a^{kj}\partial_ka^{li}-a^{kl}
\partial_ka^{ij})+\partial_ta^{ij}
+A(\frac{d+2b\partial_k a^{kl}x_l}{t}+d)\\
=&-\frac{4b}{t}A^2-\frac{2bx_l}{t}(a^{ki}\partial_ka^{lj}+a^{kj}\partial_ka^{li}-a^{kl}
\partial_ka^{ij}-a^{ij}\partial_k a^{kl})+\partial_ta^{ij}+Ad(\frac{1}{t}+1)\\
\geq&-\frac{4b}{t}\Lambda^2I_n-\frac{2bx_l}{t}(a^{ki}\partial_ka^{lj}+a^{kj}
\partial_ka^{li}-a^{kl}\partial_ka^{ij}-a^{ij}\partial_k a^{kl})+\partial_ta^{ij}+\lambda d(\frac{1}{t}+1)I_n.
\end{split}
\end{eqnarray*}
Next we estimate the lower bounds of the matrices in the right side of the above inequality.

We just need to estimate matrix $x_la^{ki}\partial_ka^{lj}$ and
$\partial_ta^{ij}$. For any $\xi\in\mathbb{R}^n$,
$$|x_la^{ki}\partial_ka^{lj}\xi_i\xi_j|\leq n^2|x|\Lambda\frac{E}{|x|}\sum_{i,j}|\xi_i||\xi_j|\leq
n^3\Lambda E|\xi|^2,$$
then
$$-n^3\Lambda EI_n\leq x_la^{ki}\partial_ka^{lj}\leq n^3\Lambda EI_n;$$
and $$|\partial_ta^{ij}\xi_i\xi_j|\leq M\sum_{i,j}|\xi_i||\xi_j|\leq Mn|\xi|^2,$$
then
$$\partial_ta^{ij}\geq-MnI_n.$$
Thus we have
\begin{eqnarray*}
2D_G+A(\frac{\partial_t G-\tilde{\Delta}G}{G}-F)
\geq[-\frac{4b}{t}\Lambda^2-\frac{8b}{t}n^3\Lambda E-Mn+\lambda d(\frac{1}{t}+1)]I_n,
\end{eqnarray*}
if we choose $d=d(n,\Lambda,\lambda,M,E)$ large enough, then
$$2D_G+A(\frac{\partial_t G-\tilde{\Delta}G}{G}-F)\geq(\frac{1}{t}+1)I_n.$$

\emph{Estimate of $\partial_tF+F(\frac{\partial_tG-\tilde{\Delta}G}{G}-F)$.}\\

By direct calculations we have
\begin{eqnarray*}
\begin{split}
&\partial_tF+F(\frac{\partial_tG-\tilde{\Delta}G}{G}-F)\\
=&\frac{(d-2+2b\partial_ia^{ij}x_j)(b|x|^2-4b^2a^{ij}x_ix_j+K)}{t^3}\\
&+\frac{(db|x|^2-4db^2a^{ij}x_ix_j-4b^2\partial_ta^{ij}x_ix_j)+Kd-(d-2ba^{ii})(d-1+2b\partial_ia^{ij}x_j)}{t^2}\\
&-\frac{d(2d-2ba^{ii}+2b\partial_ia^{ij}x_j)-2b\partial_ta^{ii}}{t}-d^2\\
&+2\gamma Kt^{-K-2}(K+1-d-td-2b\partial_ia^{ij}x_j).\\
\end{split}
\end{eqnarray*}
Noticing that
$$|\partial_ia^{ij}x_j|\leq n^2\frac{E}{|x|}|x|=n^2E,$$
$$a^{ij}x_ix_j\leq\Lambda|x|^2,~~~|\partial_ta^{ij}x_ix_j|\leq M\sum_{i,j}|x_i||x_j|\leq Mn|x|^2,$$
then we have
\begin{eqnarray*}
\begin{split}
&\partial_tF+F(\frac{\partial_tG-\tilde{\Delta}G}{G}-F)\\
\geq& \frac{(d-C)[(b-4b^2\Lambda)|x|^2+K]}{t^3}+\frac{(db-4db^2\Lambda-C)|x|^2+Kd-(d+C)^2}{t^2}\\
&-\frac{d(2d+C)+C}{t}-d^2+2\gamma Kt^{-K-2}(K+1-3d-C).
\end{split}
\end{eqnarray*}
Recall that $b=\frac{1}{8\Lambda}$, and we choose $d$ large enough,
then
\begin{eqnarray*}
\begin{split}
&\partial_tF+F(\frac{\partial_tG-\tilde{\Delta}G}{G}-F)\\
\geq& \frac{d(\frac{b}{2}|x|^2+K)}{2t^3}+\frac{(db/2-C)|x|^2+Kd-2d^2}{t^2}-\frac{3d^2}{t}-
d^2+2\gamma Kt^{-K-2}(K-4d)\\
\geq& \frac{d(b|x|^2+K)}{4t^3}+\frac{(db/2-C)|x|^2+Kd-12d^2}{t^2}+2\gamma Kt^{-K-2}(K-4d),
\end{split}
\end{eqnarray*}
Since $K=12d$, then
$$\partial_tF+F(\frac{\partial_tG-\tilde{\Delta}G}{G}-F)\geq \frac{db(|x|^2+1)}{4t^3}.$$

\emph{Estimate of $\tilde{\triangle}F_0$}.\\

In order to estimate $\tilde{\Delta} F_0$ and $|\nabla(F-F_0)|$, we
need some estimates about $\{a^{ij}_\epsilon\}$ which we put in
Appendix A.

In fact, $\{a^{ij}_\epsilon\}$ satisfy the following properties:
\begin{equation}\label{indexA}
\begin{split}
&i)~\lambda|\xi|^2\leq a^{ij}_\epsilon(x,t)\xi_i\xi_j\leq \Lambda|\xi|^2, ~\forall \xi\in\mathbb{R}^n;\\
&ii)~|\nabla a^{ij}_\epsilon(x,t)|\leq M;~~|\nabla a^{ij}_\epsilon(x,t)|\leq\frac{2E}{|x|}~when~|x|\geq1;\\
&iii)~|a^{ij}_\epsilon(x,t)-a^{ij}(x,t)|\leq 2\Lambda ;
~~|a^{ij}_\epsilon(x,t)-a^{ij}(x,t)|\leq\frac{E}{|x|}~when~|x|\geq1;\\
&iv)~|\partial_{kl}a^{ij}_\epsilon(x,t)|\leq c(n)M;~~|\partial_{kl}a^{ij}_\epsilon(x,t)|\leq
\frac{c(n)E}{|x|}~when~|x|\geq1.\\
\end{split}
\end{equation}
Direct calculations show that
\begin{eqnarray}
\begin{split}
\tilde{\triangle}F_0=&\frac{b}{t^2}\tilde{\triangle}(|x|^2)-\frac{4b^2}{t^2}
\tilde{\triangle}(a^{ij}_\epsilon x_ix_j)+\frac{2b}{t}\tilde{\triangle}(a^{ii}_\epsilon)\\
=&\frac{2b}{t^2}(a^{ii}+\partial_ia^{ij}x_j)-\frac{4b^2}{t^2}
[(a^{kl}\partial_{kl}a^{ij}_\epsilon+\partial_ka^{kl}\partial_la^{ij}_\epsilon)x_ix_j\\
&+2(\partial_ka^{kj}a^{ij}_\epsilon+2a^{kj}\partial_ka^{ij}_\epsilon)x_i+2a^{ij}a^{ij}_\epsilon]
+\frac{2b}{t}(a^{kl}\partial_{kl}a^{ii}_\epsilon+\partial_ka^{kl}\partial_la^{ii}_\epsilon).
\end{split}
\end{eqnarray}
Now we estimate the terms in the right side of the above identity.
By (\ref{indexA}) we have
$$|a^{ij}_\epsilon|,~~|\nabla a^{ij}_\epsilon|~~\text{and}~~|\nabla^2a^{ij}_\epsilon|~~\text{are all bounded},$$
then
$$|a^{ii}+\partial_ia^{ij}x_j|\leq C(1+|x|);$$
$$|(a^{kl}\partial_{kl}a^{ij}_\epsilon+\partial_ka^{kl}\partial_la^{ij}_\epsilon)x_ix_j|\leq C|x|^2;$$
$$|(\partial_ka^{kj}a^{ij}_\epsilon+2a^{kj}\partial_ka^{ij}_\epsilon)x_i|\leq C|x|;$$
$$|a^{kl}\partial_{kl}a^{ii}_\epsilon+\partial_ka^{kl}\partial_la^{ii}_\epsilon|\leq C.$$
Thus
$$
\tilde{\triangle}F_0\geq-\frac{C(|x|+1)}{t^2}-\frac{C(|x|^2+|x|+1)}{t^2}-\frac{C}{t}\geq-\frac{C(|x|^2+1)}{t^3}.
$$

\emph{Estimate of $|\nabla(F-F_0)|$.}\\

Since
$$F-F_0=\frac{4b^2(a^{ij}_\epsilon-a^{ij})x_ix_j}{t^2}-\frac{2b(a^{ii}_\epsilon-a^{ii})}{t},$$
then
\begin{eqnarray*}
\begin{split}
|\nabla(F-F_0)|=&|\frac{4b^2}{t^2}(\nabla a^{ij}_\epsilon-\nabla a^{ij})x_ix_j+
\frac{8b^2}{t^2}(a^{ij}_\epsilon-a^{ij})x_i\nabla x_j-\frac{2b}{t}(\nabla a^{ii}_\epsilon-\nabla a^{ii})|\\
\leq& \frac{4b^2}{t^2}|\nabla a^{ij}_\epsilon-\nabla a^{ij}| |x_i||x_j|
+\frac{8b^2}{t^2}\sum_j|a^{ij}_\epsilon-a^{ij}||x_i|+\frac{2b}{t}|\nabla a^{ii}_\epsilon-\nabla a^{ii}|.
\end{split}
\end{eqnarray*}
We now estimate the terms in the right side of the above inequality.\\
By $ii)$ of (\ref{indexA}), when $|x|<1$,
$$|\nabla a^{ij}_\epsilon-\nabla a^{ij}| |x_i||x_j|\leq 2M\sum_{i,j}|x_i||x_j|\leq 2Mn|x|^2\leq 2Mn,$$
and when $|x|\geq1$,
$$|\nabla a^{ij}_\epsilon-\nabla a^{ij}| |x_i||x_j|\leq (\frac{2E}{|x|}+
\frac{E}{|x|})\sum_{i,j}|x_i||x_j|\leq \frac{3E}{|x|}n|x|^2=3nE|x|.
$$
Then we have
$$|\nabla a^{ij}_\epsilon-\nabla a^{ij}| |x_i||x_j|\leq C(|x|+1).$$
By $ii)$ and $iii)$ of (\ref{indexA}) we have
$$\sum_j|a^{ij}_\epsilon-a^{ij}||x_i|\leq 2n^2\Lambda |x|,$$
$$|\nabla a^{ii}_\epsilon-\nabla a^{ii}|\leq 2Mn.$$
With the above three estimates we have
$$
|\nabla(F-F_0)|\leq \frac{C(|x|+1)}{t^2}+\frac{C|x|}{t^2}+\frac{C}{t}\leq\frac{C(|x|+1)}{t^2}.
$$
Thus we proved Lemma \ref{estimates1}.

\subsection{Proof of Proposition \ref{Prop-C2}.}

Before we prove proposition \ref{Prop-C2}, we need to prove a result as another
version of Corollary \ref{CIforgeneral}.\\

In (\ref{pc0.2}), we let $\Phi=\gamma (t^{-K}-1)x_n^{\alpha}-\frac{b\psi+K}{2t}$,
$G=e^{2\Phi}$, $v=e^{\Phi}u$ and we denote
$$B=2D_G+A(\frac{\partial_tG-\tilde{\Delta}G}{G}-F).$$
Then the third term of the left hand side of (\ref{pc0.2}) is
\begin{eqnarray*}
\begin{split}
&-\int_Qu\langle A\nabla(F-F_0),\nabla u\rangle e^{2\Phi}dxdt\\
=&-\int_Qv\langle A\nabla(F-F_0),\nabla v-\nabla \Phi v\rangle dxdt\\
=&-\int_Qv\langle A\nabla(F-F_0),\nabla v\rangle dxdt+\int_Q\langle A\nabla(F-F_0),\nabla \Phi \rangle v^2dxdt.\\
\end{split}
\end{eqnarray*}
We use the above identity and rewrite (\ref{pc0.2}) as
\begin{equation}\label{pc2.1}
\begin{split}
\frac{1}{2}\int_QM_1v^2dxdt+&\int_Q\langle B\nabla u,\nabla u\rangle e^{2\Phi}dxdt-
\int_Qv\langle A\nabla(F-F_0),\nabla v\rangle dxdt\\
&=2\int_QLu(Pu-Lu)e^{2\Phi}dxdt,
\end{split}
\end{equation}
where
$$M_1=\partial_tF+F(\frac{\partial_tG-\tilde{\Delta}G}{G}-F)+\tilde{\Delta}F_0,$$
$$B=4AD^2\Phi A+2\partial_l\Phi(a^{ki}\partial_ka^{lj}+a^{kj}\partial_ka^{li}-a^{kl}\partial_ka^{ij})
+\partial_ta^{ij}+A(\frac{\partial_tG-\tilde{\Delta}G}{G}-F).$$
We rewrite $\Phi$ as the following:
\begin{eqnarray*}
\begin{split}
&\Phi=\Phi_1+\Phi_2,\\
&\Phi_1=\gamma f(t)x_n^\alpha, ~~f(t)=t^{-K}-1,\\
&\Phi_2=-\frac{b\psi+K}{2t}.
\end{split}
\end{eqnarray*}
The function $\psi$ has the following properties which we will prove
in Appendix B:
\begin{equation}\label{pc2.2}
\begin{split}
&i)~\psi\geq\frac{|x|^2}{2};\\
&ii)~D^2\psi\leq C(\frac{\Lambda}{\lambda})I_n;\\
&iii)~|\nabla \psi|\leq 4(\frac{\Lambda}{\lambda}+1)^2|x|, ~|\nabla^k\psi|\leq
\frac{C(n,\frac{\Lambda}{\lambda})}{|x|^{k-2}},\quad k=2,3,4;\\
&iv)~a^{ni}\partial_i\psi\leq C(\Lambda,\lambda)x_n.
\end{split}
\end{equation}
By direct calculations we have
$$\frac{\partial_tG-\tilde{\Delta}G}{G}=2\partial_t\Phi-2a^{ij}\partial_{ij}\Phi-
2\partial_ia^{ij}\partial_j\Phi-4\langle A\nabla\Phi,\nabla\Phi\rangle.$$
Let
\begin{equation}\label{choiceofF}
F=2\partial_t\Phi-2a^{ij}\partial_{ij}\Phi-4\langle A\nabla\Phi,\nabla\Phi\rangle-H,
\end{equation}
where $H$ is a positive smooth function to be determined. Let
$$F_0=2\partial_t\Phi-2a^{ij}_\epsilon\partial_{ij}\Phi-4a^{ij}_\epsilon\partial_i\Phi\partial_j\Phi-H.$$
We estimate matrix $B$ first. Direct calculations show that
\begin{eqnarray*}
\begin{split}
B&=4AD^2\Phi A+2\partial_l\Phi(a^{ki}\partial_ka^{lj}+a^{kj}\partial_ka^{li}-a^{kl}\partial_ka^{ij})
+\partial_ta^{ij}+A(-2\partial_ka^{kl}\partial_l\Phi+H)\\
&=4AD^2\Phi A+2\partial_l\Phi(a^{ki}\partial_ka^{lj}+a^{kj}\partial_ka^{li}-a^{kl}\partial_ka^{ij}-a^{ij}\partial_ka^{kl})
+\partial_ta^{ij}+HA.
\end{split}
\end{eqnarray*}
We estimate the lower bounds of the matrices in the right side of the above identity.\\
First, by $ii)$ of (\ref{pc2.2}) we have
$$
AD^2\Phi A=AD^2\Phi_1 A-\frac{b}{2t}AD^2\psi A
\geq\partial_{nn}\Phi_1a^{in}a^{nj}-\frac{C}{t}I_n.
$$
Second, we estimate matrix $\partial_l\Phi a^{ki}\partial_ka^{lj}$ and $\partial_ta^{ij}$.\\
For any $\xi\in\mathbb{R}^n$,
\begin{eqnarray*}
\begin{split}
|\partial_l\Phi a^{ki}\partial_ka^{lj}\xi_i\xi_j|\leq& n\Lambda\frac{E}{|x|}
\sum_{l,i,j}|\partial_l\Phi|\xi_i||\xi_j|\\
\leq& \frac{n^2\Lambda E}{|x|}|\xi|^2\sum_l|\partial_l\Phi|\\
\leq& \frac{n^2\Lambda E}{|x|}(\partial_n\Phi_1+\frac{C}{t}|\nabla\psi|)|\xi|^2,
\end{split}
\end{eqnarray*}
and by $iii)$ of (\ref{pc2.2}), $|\nabla \psi|\leq C|x|$, then
$$|\partial_l\Phi a^{ki}\partial_ka^{lj}\xi_i\xi_j|\leq(n^2\Lambda E
\frac{\partial_n\Phi_1}{|x|}+\frac{C}{t})|\xi|^2,$$
thus
$$-(n^2\Lambda E\frac{\partial_n\Phi_1}{|x|}+\frac{C}{t})I_n\leq
\partial_l\Phi a^{ki}\partial_ka^{lj}\leq (n^2\Lambda E\frac{\partial_n\Phi_1}{|x|}+\frac{C}{t})I_n;$$
again,
$$|\partial_ta^{ij}\xi_i\xi_j|\leq M\sum_{i,j}|\xi_i||\xi_j|\leq Mn|\xi|^2,$$
then
$$\partial_ta^{ij}\geq-MnI_n.$$
Consequently,
\begin{equation}\label{B-tildeB}
\begin{split}
B\geq&4\partial_{nn}\Phi_1a^{in}a^{nj}-\frac{C}{t}I_n-8(n^2\Lambda E\frac{\partial_n\Phi_1}{|x|}+
\frac{C}{t})I_n-MnI_n+HA\\
\geq&4\partial_{nn}\Phi_1a^{in}a^{nj}-(8n^2\Lambda E\frac{\partial_n\Phi_1}{|x|}+
\frac{C}{t})I_n+HA+\frac{1}{t}I_n\\
\equiv&\tilde{B}+\frac{1}{t}I_n,\\
\end{split}
\end{equation}
where
\begin{equation}\label{tildeB}
\tilde{B}=4\partial_{nn}\Phi_1a^{in}a^{nj}-(8n^2\Lambda E\frac{\partial_n\Phi_1}{|x|}+\frac{C}{t})I_n+HA.
\end{equation}
To make $\tilde{B}$ positive, we choose
\begin{equation}\label{choiceofH}
H=16n^2\frac{\Lambda}{\lambda}E\frac{\partial_n\Phi_1}{|x|}+\frac{d}{t},
\end{equation}
where $d$ is a positive constant to be determined.\\
Since $\tilde{B}$ is differentiable, then by (\ref{B-tildeB}) we have
\begin{equation}\label{pc2.3}
\begin{split}
&\int_Q\langle B\nabla u,\nabla u\rangle e^{2\Phi}dxdt\\
\geq&\int_Q\langle \tilde{B}\nabla u,\nabla u\rangle e^{2\Phi}dxdt+\int_Qe^{2\Phi}\frac{|\nabla u|^2}{t}dxdt\\
=&\int_Q\langle \tilde{B}\nabla v,\nabla v\rangle dxdt+\int_Q[\langle\tilde{B}\nabla \Phi,\nabla \Phi\rangle
+div(\tilde{B}\nabla\Phi)]v^2dxdt+\int_Qe^{2\Phi}\frac{|\nabla u|^2}{t}dxdt.\\
\end{split}
\end{equation}
By (\ref{pc2.1}), (\ref{pc2.3}) and the Cauchy inequality, we have
\begin{equation}\label{pc2.4}
\begin{split}
&\int_Qe^{2\Phi}\frac{|\nabla u|^2}{t}dxdt+\int_Q\langle \tilde{B}\nabla v,\nabla v\rangle dxdt+\int_QM_2v^2dxdt\\
&-\int_Qv\langle A\nabla(F-F_0),\nabla v\rangle dxdt\leq \int_Qe^{2\Phi}|Pu|^2dxdt,
\end{split}
\end{equation}
where
$$M_2=\langle\tilde{B}\nabla \Phi,\nabla \Phi\rangle+div(\tilde{B}\nabla\Phi)
+\frac{1}{2}\partial_tF+\frac{1}{2}F(\frac{\partial_tG-\tilde{\Delta}G}{G}-F)+\frac{1}{2}\tilde{\Delta}F_0.$$

We use inequality (\ref{pc2.4}) to prove Proposition \ref{Prop-C2}.
We also need some estimates which we list in the following lemma. We will prove this lemma later.

\begin{lemma}\label{estimates2}
Set $b=\frac{1}{64\Lambda(\frac{\Lambda}{\lambda}+1)^4}$, $E_0=
\frac{\lambda}{16n^2\frac{\Lambda}{\lambda}(\frac{\Lambda}{\lambda}+1)}$,
$\alpha=1+\frac{E}{E_0}$ and $K=13\frac{\Lambda}{\lambda}d$. We take $d=d(n,\Lambda,\lambda,M,E)$
large enough, when $E<E_0$, we have
\begin{equation}\label{E5}
\tilde{B}\geq8n^2\Lambda E(\frac{\partial_n\Phi_1}{|x|}+\frac{1}{t})I_n;
\end{equation}
\begin{equation}\label{E6}
M_2\geq2[(\alpha-1)\lambda^2-(16n^2\frac{\Lambda}{\lambda}+8n^2+4n)\Lambda E]\frac{(\partial_n\Phi_1)^3}{|x|}
+\frac{bd|x|^2}{16t^3}+\frac{1}{t^3};
\end{equation}
\begin{equation}\label{E7}
|\nabla(F-F_0)|\leq32nE[\frac{(\partial_n\Phi_1)^2}{|x|}+\frac{C|x|}{t^2}].
\end{equation}
\end{lemma}

Then by applying (\ref{E5}) in Lemma 3.4,  we have
\begin{equation}\label{pc2.5}
\int_Q\langle \tilde{B}\nabla v,\nabla v\rangle dxdt\geq8n^2\Lambda E\int_Q(\frac{\partial_n\Phi_1}{|x|}
+\frac{1}{t})|\nabla v|^2dxdt.
\end{equation}
By (\ref{E6}) we have
\begin{equation}\label{pc2.15}
\begin{split}
\int_QM_2v^2dxdt\geq&2\int_Q[(\alpha-1)\lambda^2-(16n^2\frac{\Lambda}{\lambda}+8n^2+4n)\Lambda E]
\frac{(\partial_n\Phi_1)^3}{|x|}v^2dxdt\\
&+\int_Q(\frac{bd|x|^2}{16t^3}+\frac{1}{t^3})v^2dxdt.
\end{split}
\end{equation}
By (\ref{E7})
\begin{eqnarray*}
\begin{split}
|\int_Qv\langle A\nabla(F-F_0),\nabla v\rangle dxdt|\leq&\Lambda\int_Q|\nabla(F-F_0)||v||\nabla v|dxdt\\
\leq&32n\Lambda
E\int_Q[\frac{(\partial_n\Phi_1)^2}{|x|}+\frac{C|x|}{t^2}]|v||\nabla
v|dxdt.
\end{split}
\end{eqnarray*}
Using the Cauchy inequality, we have
\begin{equation}\label{pc2.17}
\begin{split}
|\int_Qv\langle A\nabla(F-F_0),\nabla v\rangle dxdt|
\leq&32\Lambda E\int_Q[\frac{(\partial_n\Phi_1)^3}{|x|}+\frac{C|x|^2}{t^3}]v^2dxdt\\
&+8n^2\Lambda E\int_Q(\frac{\partial_n\Phi_1}{|x|}+\frac{1}{t})|\nabla v|^2dxdt.
\end{split}
\end{equation}
Because of (\ref{pc2.4}), (\ref{pc2.5}), (\ref{pc2.15}) and (\ref{pc2.17}), we have
\begin{eqnarray*}
\begin{split}
\int_Qe^{2\Phi}|Pu|^2dxdt\geq&\int_Q2[(\alpha-1)\lambda^2-(16n^2\frac{\Lambda}{\lambda}+8n^2+4n+16)
\Lambda E]\frac{(\partial_n\Phi_1)^3}{|x|}v^2dxdt\\
&+\int_Q[\frac{(bd-C)|x|^2}{16t^3}+\frac{1}{t^3}]v^2dxdt+\int_Qe^{2\Phi}\frac{|\nabla u|^2}{t}dxdt,
\end{split}
\end{eqnarray*}
Since $$16n^2\frac{\Lambda}{\lambda}+8n^2+4n+16\leq16n^2(\frac{\Lambda}{\lambda}+1),$$
and if we take $d$ large enough, then
\begin{eqnarray*}
\begin{split}
\int_Qe^{2\Phi}|Pu|^2dxdt\geq&\int_Q2[(\alpha-1)\lambda^2-16n^2(\frac{\Lambda}{\lambda}+1)
\Lambda E]\frac{(\partial_n\Phi_1)^3}{|x|}v^2dxdt\\
&+\int_Qe^{2\Phi}(\frac{|u|^2}{t^3}+\frac{|\nabla u|^2}{t})dxdt.
\end{split}
\end{eqnarray*}
Notice that
$$(\alpha-1)\lambda^2-16n^2(\frac{\Lambda}{\lambda}+1)\Lambda E=0,$$
thus
$$\int_Qe^{2\Phi}|Pu|^2dxdt\geq\int_Qe^{2\Phi}(\frac{|u|^2}{t^3}+\frac{|\nabla u|^2}{t})dxdt.$$
Thus we proved Carleman inequality (\ref{C2}).\\

\emph{Proof of Lemma \ref{estimates2}}. \\

\emph{Estimate of $\tilde{B}$.}\\

By (\ref{tildeB}) and (\ref{choiceofH}) we have
$$
\tilde{B}\geq[-(8n^2\Lambda
E\frac{\partial_n\Phi_1}{|x|}+\frac{C}{t})+\lambda H]I_n
=(8n^2\Lambda E\frac{\partial_n\Phi_1}{|x|}+\frac{\lambda
d-C}{t})I_n,
$$
if we take $d=d(n,\Lambda,\lambda,M,E)$ large enough, then
$$\tilde{B}\geq8n^2\Lambda E(\frac{\partial_n\Phi_1}{|x|}+\frac{1}{t})I_n.$$

\emph{Estimate of $M_2$}.\\

In order to estimate $M_2$, we have to divide $M_2$ into several parts and estimate each of them:
$$M_2=J_1+J_2+J_3+J_4+J_5+J_6,$$
where
\begin{equation}\label{divideJ}
\begin{split}
J_1=&4\partial_{nn}\Phi_1(a^{ni}\partial_i\Phi)^2;\\
J_2=&-(8n^2\Lambda E\frac{\partial_n\Phi_1}{|x|}+\frac{C}{t})|\nabla\Phi|^2-(H-4\partial_ia^{ij}\partial_j\Phi)
\langle A\nabla\Phi,\nabla\Phi\rangle\\
&+\langle A\nabla H,\nabla\Phi\rangle-8n^2\Lambda E\langle\nabla(\frac{\partial_n\Phi_1}{|x|}),\nabla\Phi\rangle
-(8n^2\Lambda E\frac{\partial_n\Phi_1}{|x|}+\frac{C}{t})\Delta\Phi;\\
J_3=&4\partial_{nn}\Phi_1\partial_i(a^{in}a^{nj})\partial_j\Phi+4\partial_n^3\Phi_1a^{nn}a^{nj}\partial_j\Phi
+4\partial_{nn}\Phi_1a^{in}a^{nj}\partial_{ij}\Phi;\\
J_4=&\partial_{tt}\Phi+\partial_t\Phi(H-2\partial_ia^{ij}\partial_j\Phi)-
\partial_ta^{ij}\partial_{ij}\Phi-a^{ij}\partial_{ijt}\Phi\\
&+2\partial_ia^{ij}\partial_j\Phi(H+a^{ij}\partial_{ij}\Phi)-\frac{1}{2}\partial_t H-\frac{1}{2}H^2;\\
J_5=&-2\partial_t\langle A\nabla\Phi,\nabla\Phi\rangle;\\
J_6=&\frac{1}{2}\tilde{\Delta}F_0.
\end{split}
\end{equation}

\emph{Estimate of $J_1$}.\\
$$J_1=4\partial_{nn}\Phi_1(a^{nn}\partial_n\Phi_1+a^{ni}\partial_i\Phi_2)^2,$$
By the Cauchy inequality $(a+b)^2\geq\delta a^2-\frac{\delta}{1-\delta}b^2,~0<\delta<1$, we have
\begin{eqnarray*}
\begin{split}
J_1&\geq4\partial_{nn}\Phi_1[\frac{x_n}{2|x|}(a^{nn}\partial_n\Phi_1)^2-
\frac{x_n}{2|x|-x_n}(a^{ni}\partial_i\Phi_2)^2]\\
&\geq4\partial_{nn}\Phi_1[\frac{x_n}{2|x|}\lambda^2(\partial_n\Phi_1)^2-
\frac{x_n}{|x|}(a^{ni}\partial_i\Phi_2)^2]\\
&=2(\alpha-1)\frac{\partial_n\Phi_1}{|x|}[\lambda^2(\partial_n\Phi_1)^2-2(a^{ni}\partial_i\Phi_2)^2].
\end{split}
\end{eqnarray*}
Again by the Cauchy inequality we have
$$(a^{ni}\partial_i\Phi_2)^2\leq\sum_i(a^{ni})^2|\nabla\Phi_2|^2\leq \frac{\Lambda^2b^2}{4t^2}|\nabla\psi|^2,$$
then
\begin{eqnarray*}
\begin{split}
J_1\geq&2(\alpha-1)\frac{\partial_n\Phi_1}{|x|}[\lambda^2(\partial_n\Phi_1)^2-
\frac{\Lambda^2b^2}{2t^2}|\nabla\psi|^2]\\
=&2(\alpha-1)\lambda^2\frac{(\partial_{n}\Phi_1)^3}{|x|}-\frac{C}{t^2}
\frac{\partial_{n}\Phi_1}{|x|}|\nabla\psi|^2.
\end{split}
\end{eqnarray*}
By $iii)$ of (\ref{pc2.2}), $|\nabla \psi|\leq C|x|$, then
$$J_1\geq2(\alpha-1)\lambda^2\frac{(\partial_{n}\Phi_1)^3}{|x|}-\frac{C}{t^2}\partial_{n}\Phi_1|x|.$$
Using the Cauchy inequality, we obtain
\begin{equation}\label{pc2.6}
J_1\geq2(\alpha-1)\lambda^2\frac{(\partial_{n}\Phi_1)^3}{|x|}-
\frac{C}{t}(\partial_{n}\Phi_1)^2-\frac{C}{t^3}|x|^2.\\
\end{equation}

\emph{Estimate of $J_2$}.\\

Let's estimate $\partial_ia^{ij}\partial_j\Phi$ first. We will use
 this estimate afterwards again.
$$|\partial_ia^{ij}\partial_j\Phi|\leq\frac{nE}{|x|}\sum_j|\partial_j\Phi|\leq
\frac{nE}{|x|}(\partial_n\Phi_1+\frac{C}{t}|\nabla\psi|).$$ Recall
that $|\nabla\psi|\leq C|x|$, then
\begin{equation}\label{ofenuse}
|\partial_ia^{ij}\partial_j\Phi|\leq nE\frac{\partial_n\Phi_1}{|x|}+\frac{C}{t}.
\end{equation}
By (\ref{divideJ}) we have
\begin{eqnarray*}
\begin{split}
J_2\geq&-(8n^2\Lambda E\frac{\partial_n\Phi_1}{|x|}+\frac{C}{t})|\nabla\Phi|^2-
[(16n^2\frac{\Lambda}{\lambda}+4n)E\frac{\partial_n\Phi_1}{|x|}+\frac{d+C}{t}]
\Lambda|\nabla\Phi|^2\\
&-\Lambda|\nabla H||\nabla\Phi|-8n^2\Lambda E|\nabla(\frac{\partial_n\Phi_1}{|x|})||\nabla\Phi|
-(8n^2\Lambda E\frac{\partial_n\Phi_1}{|x|}+\frac{C}{t})\Delta\Phi\\
=&-[(16n^2\frac{\Lambda}{\lambda}+8n^2+4n)\Lambda E\frac{\partial_n\Phi_1}{|x|}+
\frac{d\Lambda+C}{t}]|\nabla\Phi|^2\\
&-(16n^2\frac{\Lambda}{\lambda}+8n^2)\Lambda E|\nabla(\frac{\partial_n\Phi_1}{|x|})||\nabla\Phi|-
(8n^2\Lambda E\frac{\partial_n\Phi_1}{|x|}+\frac{C}{t})\Delta\Phi.
\end{split}
\end{eqnarray*}
Because
$$|\nabla(\frac{\partial_n\Phi_1}{|x|})|\leq\frac{\partial_n\Phi_1}{x_n|x|}\leq \partial_n\Phi_1,$$
and we have by $iii)$ of (\ref{pc2.2}) that
$$\Delta\Phi=\partial_{nn}\Phi_1-\frac{b}{2t}\Delta\psi\leq\partial_{nn}\Phi_1+\frac{C}{t},$$
then
\begin{eqnarray*}
\begin{split}
J_2\geq&-[(16n^2\frac{\Lambda}{\lambda}+8n^2+4n)\Lambda E\frac{\partial_n\Phi_1}{|x|}+
\frac{d\Lambda+C}{t}]|\nabla\Phi|^2
-C\partial_n\Phi_1|\nabla\Phi|\\
&-(8n^2\Lambda E\frac{\partial_n\Phi_1}{|x|}+\frac{C}{t})(\partial_{nn}\Phi_1+\frac{C}{t})\\
\geq&-[(16n^2\frac{\Lambda}{\lambda}+8n^2+4n)\Lambda E\frac{\partial_n\Phi_1}{|x|}+
\frac{d\Lambda+C}{t}]|\nabla\Phi|^2-C\partial_n\Phi_1|\nabla\Phi|\\
&-C(\partial_n\Phi_1)^2-\frac{C}{t}\partial_n\Phi_1-\frac{C}{t^2}.
\end{split}
\end{eqnarray*}
By the Cauchy inequality we have
$$C\partial_n\Phi_1|\nabla\Phi|\leq C(\partial_n\Phi_1)^2+C|\nabla\Phi|^2,
~~\frac{C}{t}\partial_n\Phi_1\leq C(\partial_n\Phi_1)^2+\frac{C}{t^2},$$
then
\begin{eqnarray*}
\begin{split}
J_2\geq&-[(16n^2\frac{\Lambda}{\lambda}+8n^2+4n)\Lambda E\frac{\partial_n\Phi_1}{|x|}+
\frac{d\Lambda+C}{t}]|\nabla\Phi|^2-C(\partial_{n}\Phi_1)^2-\frac{C}{t^2}\\
\geq&-2[(16n^2\frac{\Lambda}{\lambda}+8n^2+4n)\Lambda E\frac{\partial_n\Phi_1}{|x|}+
\frac{d\Lambda+C}{t}][(\partial_n\Phi_1)^2+\frac{b^2}{4t^2}|\nabla\psi|^2]
-C(\partial_{n}\Phi_1)^2-\frac{C}{t^2}\\
\geq&-2(16n^2\frac{\Lambda}{\lambda}+8n^2+4n)\Lambda E\frac{(\partial_n\Phi_1)^3}{|x|}
-\frac{2d\Lambda+C}{t}(\partial_{n}\Phi_1)^2-\frac{C}{t^2}\frac{\partial_n\Phi_1}{|x|}|\nabla\psi|^2\\
&-\frac{(d\Lambda+C)b^2}{2t^3}|\nabla\psi|^2-\frac{C}{t^2}.
\end{split}
\end{eqnarray*}
Taking in account that $|\nabla\psi|\leq C|x|$ and using the Cauchy inequality, we have that
$$\frac{C}{t^2}\frac{\partial_n\Phi_1}{|x|}|\nabla\psi|^2\leq\frac{C}{t^2}
\partial_n\Phi_1|x|\leq \frac{C}{t}(\partial_n\Phi_1)^2+\frac{C}{t^3}|x|^2,$$
then
\begin{equation}\label{pc2.7}
\begin{split}
J_2\geq&-2(16n^2\frac{\Lambda}{\lambda}+8n^2+4n)\Lambda E\frac{(\partial_n\Phi_1)^3}{|x|}-
\frac{2d\Lambda+C}{t}(\partial_n\Phi_1)^2
-\frac{C}{t^3}|x|^2\\
&-\frac{db^2\Lambda+C}{2t^3}|\nabla\psi|^2-\frac{C}{t^2}.\\
\end{split}
\end{equation}
In the following, we always use the fact that
$$|\nabla^k\psi|\leq \frac{C}{|x|^{k-2}},~~k=1,2,3,4.$$

\emph{Estimate of $J_3$}.\\

$$J_3=4\partial_{nn}\Phi_1(\partial_ia^{in}a^{nj}+a^{in}\partial_ia^{nj})\partial_j\Phi+
4\partial_n^3\Phi_1a^{nn}a^{nj}\partial_j\Phi
+4\partial_{nn}\Phi_1a^{in}a^{nj}\partial_{ij}\Phi.$$
Next we estimate the terms of $J_3$.
$$|(\partial_ia^{in}a^{nj}+a^{in}\partial_ia^{nj})\partial_j\Phi|
\leq C|\nabla\Phi|\leq C(\partial_n\Phi_1+\frac{|\nabla\psi|}{t})\leq C(\partial_n\Phi_1+\frac{|x|}{t});$$
$$|a^{nn}a^{nj}\partial_j\Phi|\leq C|\nabla\Phi|\leq C(\partial_n\Phi_1+
\frac{|\nabla\psi|}{t})\leq C(\partial_n\Phi_1+\frac{|x|}{t});$$
$$|a^{in}a^{nj}\partial_{ij}\Phi|\leq C|\nabla^2\Phi|\leq C(\partial_{nn}\Phi_1+
\frac{|\nabla^2\psi|}{t})\leq C(\partial_{nn}\Phi_1+\frac{1}{t}).$$
Combining the above estimates, we have
\begin{eqnarray*}
\begin{split}
J_3\geq&-C\partial_{nn}\Phi_1(\partial_n\Phi_1+\frac{|x|}{t})-C|\partial_n^3\Phi_1|(\partial_n\Phi_1+\frac{|x|}{t})-C\partial_{nn}\Phi_1(\partial_{nn}\Phi_1+\frac{1}{t})\\
\geq&-C[(\partial_n\Phi_1)^2+\frac{1}{t}\partial_{n}\Phi_1+\partial_n\Phi_1\frac{|x|}{t}].
\end{split}
\end{eqnarray*}
Using the Cauchy inequality, we have
\begin{equation}\label{pc2.8}
J_3\geq-C[(\partial_n\Phi_1)^2+\frac{1}{t^2}+\frac{|x|^2}{t^2}].
\end{equation}

\emph{Estimate of $J_4$}.\\
\begin{eqnarray*}
\begin{split}
J_4=&\partial_{tt}\Phi+\partial_t\Phi(H-2\partial_ia^{ij}\partial_j\Phi)-\partial_ta^{ij}
\partial_{ij}\Phi-a^{ij}\partial_{ijt}\Phi\\
&+2\partial_ia^{ij}\partial_j\Phi(H+a^{ij}\partial_{ij}\Phi)-\frac{1}{2}\partial_t H-\frac{1}{2}H^2.
\end{split}
\end{eqnarray*}
We estimate the terms of $J_4$.\\
In fact
$$\partial_{tt}\Phi=\partial_{tt}\Phi_1+\partial_{tt}\Phi_2=\partial_{tt}\Phi_1-\frac{b\psi+K}{t^3},$$

$$\partial_t\Phi(H-2\partial_ia^{ij}\partial_j\Phi)=\partial_t\Phi_1(H-2\partial_ia^{ij}\partial_j\Phi)+
\frac{b\psi+K}{2t^2}(H-2\partial_ia^{ij}\partial_j\Phi).$$
Recall (\ref{ofenuse}) $$|\partial_ia^{ij}\partial_j\Phi|\leq
nE\frac{\partial_n\Phi_1}{|x|}+\frac{C}{t},$$ and notice that
$\partial_t\Phi_1<0$, then we have
\begin{eqnarray*}
\begin{split}
&\partial_t\Phi(H-2\partial_ia^{ij}\partial_j\Phi)\\
\geq&\partial_t\Phi_1[(16n^2\frac{\Lambda}{\lambda}+2n)E\frac{\partial_n\Phi_1}{|x|}+
\frac{d+C}{t}]+\frac{b\psi+K}{2t^2}[(16n^2\frac{\Lambda}{\lambda}-2n)E\frac{\partial_n\Phi_1}{|x|}+\frac{d-C}{t}]\\
\geq&(16n^2\frac{\Lambda}{\lambda}+2n)E\frac{\partial_t\Phi_1\partial_n\Phi_1}{|x|}+
\frac{d+C}{t}\partial_t\Phi_1+\frac{(d-C)(b\psi+K)}{2t^3}.
\end{split}
\end{eqnarray*}
Because
$\frac{\partial_t\Phi_1}{|x|}\geq\frac{f'}{f}\partial_n\Phi_1$, and
if we choose $d$ large enough, then
$$\partial_t\Phi(H-2\partial_ia^{ij}\partial_j\Phi)\geq (16n^2\frac{\Lambda}{\lambda}+2n)E
\frac{f'}{f}(\partial_n\Phi_1)^2+\frac{2d}{t}\partial_t\Phi_1
+(\frac{d}{2}-C)\frac{b\psi+K}{t^3}.$$
And
$$|\partial_ta^{ij}\partial_{ij}\Phi|\leq C|\nabla^2\Phi|\leq C(\partial_{nn}\Phi_1+
\frac{|\nabla^2\psi|}{t})\leq C(\partial_n\Phi_1+\frac{1}{t}).$$
Also
$$
-a^{ij}\partial_{ijt}\Phi=-a^{nn}\partial_{nnt}\Phi_1-\frac{ba^{ij}}{2t^2}|
\partial_{ij}\psi|\geq-a^{nn}\partial_{nnt}\Phi_1-\frac{C}{t^2}|\nabla^2\psi|,
$$
because $\partial_{nnt}\Phi_1<0$, $|\nabla^2\psi|\leq C$, then
$$-a^{ij}\partial_{ijt}\Phi\geq-\frac{C}{t^2}.$$
 By (\ref{ofenuse}) we have
\begin{eqnarray*}
\begin{split}
&|2\partial_ia^{ij}\partial_j\Phi(H+a^{ij}\partial_{ij}\Phi)|\\
\leq&2(nE\frac{\partial_n\Phi_1}{|x|}+\frac{C}{t})(16n^2\frac{\Lambda}{\lambda}E
\frac{\partial_n\Phi_1}{|x|}+\frac{d}{t}+a^{nn}\partial_{nn}\Phi_1+\frac{C}{t}|\nabla^2\psi|).\\
\end{split}
\end{eqnarray*}
Notice that $|\nabla^2\psi|\leq C$, and if we choose $d$ large
enough, then
\begin{eqnarray*}
\begin{split}
|2\partial_ia^{ij}\partial_j\Phi(H+a^{ij}\partial_{ij}\Phi)|\leq&C(\partial_n\Phi_1+
\frac{1}{t})(\partial_n\Phi_1+\frac{d}{t})\\
\leq&C[(\partial_n\Phi_1)^2+\frac{d}{t}\partial_n\Phi_1+\frac{d}{t^2}].
\end{split}
\end{eqnarray*}
Using the Cauchy inequality, we obtain
$$|2\partial_ia^{ij}\partial_j\Phi(H+a^{ij}\partial_{ij}\Phi)|\leq C(\partial_n\Phi_1)^2+\frac{d^2}{t^2}.$$
And
\begin{eqnarray*}
\begin{split}
-\frac{1}{2}\partial_t H-\frac{1}{2}H^2
=&-8n^2\frac{\Lambda}{\lambda}E\frac{\partial_{nt}\Phi_1}{|x|}+\frac{d}{2t^2}-
\frac{1}{2}(16n^2\frac{\Lambda}{\lambda}E\frac{\partial_n\Phi_1}{|x|}+\frac{d}{t})^2\\
\geq&-8n^2\frac{\Lambda}{\lambda}E\frac{\partial_{nt}\Phi_1}{|x|}+\frac{d}{2t^2}-
(16n^2\frac{\Lambda}{\lambda}E\frac{\partial_n\Phi_1}{|x|})^2-\frac{d^2}{t^2}\\
\geq&-8n^2\frac{\Lambda}{\lambda}E\frac{\partial_{nt}\Phi_1}{|x|}-C(\partial_n\Phi_1)^2-\frac{d^2}{t^2}.
\end{split}
\end{eqnarray*}
Take in account that $\partial_{nt}\Phi_1<0$, then
$$-\frac{1}{2}\partial_t H-\frac{1}{2}H^2\geq-C(\partial_n\Phi_1)^2-\frac{d^2}{t^2}.$$
Combining them together, we obtain
\begin{equation}\label{pc2.12}
\begin{split}
J_4\geq&[(16n^2\frac{\Lambda}{\lambda}+2n)E\frac{f'}{f}-C](\partial_n\Phi_1)^2+\partial_{tt}\Phi_1+
\frac{2d}{t}\partial_t\Phi_1-C\partial_n\Phi_1\\
&+\frac{b\psi+K}{t^3}(\frac{d}{2}-C)-\frac{2d^2+C}{t^2}\\
\geq&[(16n^2\frac{\Lambda}{\lambda}+2n)E\frac{f'}{f}-C](\partial_n\Phi_1)^2+\partial_{tt}\Phi_1+
\frac{2d}{t}\partial_t\Phi_1-C\partial_n\Phi_1\\
&+\frac{b\psi}{t^3}(\frac{d}{2}-C)+\frac{Kd}{4t^3}-\frac{2d^2+C}{t^2}.\\
\end{split}
\end{equation}

\emph{Estimate of $J_5$}.\\
\begin{eqnarray*}
\begin{split}
J_5=&-2\partial_t\langle A\nabla\Phi_1,\nabla\Phi_1\rangle-2\partial_t\langle A\nabla\Phi_2,\nabla\Phi_2\rangle
-4\partial_t\langle A\nabla\Phi_1,\nabla\Phi_2\rangle\\
=&-2\partial_t[a^{nn}(\partial_n\Phi_1)^2]-\frac{b^2}{2}\partial_t(\frac{a^{ij}\partial_i\psi\partial_j\psi}{t^2})
+2b\partial_t(\frac{\partial_n\Phi_1}{t}a^{ni}\partial_i\psi)\\
=&-2\partial_ta^{nn}(\partial_n\Phi_1)^2-4a^{nn}\frac{f'}{f}(\partial_n\Phi_1)^2
-\frac{b^2\partial_ta^{ij}\partial_i\psi\partial_j\psi}{2t^2}+\frac{b^2a^{ij}\partial_i\psi\partial_j\psi}{t^3}\\
&+2b\frac{\partial_n\Phi_1}{t}\partial_t a^{ni}\partial_i\psi
+2b\partial_t(\frac{\partial_n\Phi_1}{t})a^{ni}\partial_i\psi.
\end{split}
\end{eqnarray*}
Because $f'<0$, $a^{nn}\geq \lambda$, $|\partial_ta^{nn}|\leq M$, and
$$|\partial_ta^{ij}\partial_i\psi\partial_j\psi|\leq M\sum_{i,j}|\partial_i\psi||\partial_j\psi|\leq
Mn|\nabla\psi|^2,$$
$$a^{ij}\partial_i\psi\partial_j\psi\geq\lambda|\nabla\psi|^2\geq0,$$
$$|\partial_t a^{ni}\partial_i\psi|\leq M\sum_i|\partial_i\psi|\leq M\sqrt{n}|\nabla\psi|,$$
then
\begin{eqnarray*}
\begin{split}
J_5\geq&-2M(\partial_n\Phi_1)^2-4\lambda\frac{f'}{f}(\partial_n\Phi_1)^2-\frac{b^2Mn|\nabla\psi|^2}{2t^2}\\
&-\frac{2b}{t}M\sqrt{n}\partial_n\Phi_1|\nabla\psi|+2\alpha b(\frac{f'}{t}-
\frac{f}{t^2})\gamma x_n^{\alpha-1}a^{ni}\partial_i\psi.
\end{split}
\end{eqnarray*}
By the Cauchy inequality we have
$$\frac{2b}{t}M\sqrt{n}\partial_n\Phi_1|\nabla\psi|\leq2M(\partial_n\Phi_1)^2+\frac{b^2Mn|\nabla\psi|^2}{2t^2},$$
then
$$J_5\geq-4M(\partial_n\Phi_1)^2-4\lambda\frac{f'}{f}(\partial_n\Phi_1)^2-
\frac{b^2Mn|\nabla\psi|^2}{t^2}+2\alpha b(\frac{f'}{t}-\frac{f}{t^2})\gamma x_n^{\alpha-1}a^{ni}\partial_i\psi.$$
By (\ref{pc2.2}), $|\nabla\psi|\leq C|x|$, $a^{ni}\partial_i\psi\leq Cx_n$, and notice that $f'<0$, then we have
\begin{equation}\label{pc2.13}
\begin{split}
J_5\geq&-4M(\partial_n\Phi_1)^2-4\lambda\frac{f'}{f}(\partial_n\Phi_1)^2-\frac{C|x|^2}{t^2}
+C(\frac{f'}{t}-\frac{f}{t^2})\gamma x_n^{\alpha}\\
=&(-4\lambda\frac{f'}{f}-4M)(\partial_n\Phi_1)^2-\frac{C|x|^2}{t^2}+C(\frac{f'}{tf}-\frac{1}{t^2})\Phi_1.\\
\end{split}
\end{equation}

\emph{Estimate of $J_6$}.\\

Recall that $$F_0=2\partial_t\Phi-2a^{ij}_\epsilon\partial_{ij}\Phi-4a^{ij}_\epsilon\partial_i\Phi\partial_j\Phi
-16n^2\frac{\Lambda}{\lambda}E\frac{\partial_n\Phi_1}{|x|}-\frac{d}{t}.$$
Direct calculations show that
\begin{eqnarray*}
\begin{split}
J_6=&\tilde{\Delta}(\partial_t\Phi)-\tilde{\Delta}[a^{ij}_\epsilon(\partial_{ij}\Phi+
2\partial_i\Phi\partial_j\Phi)]
-8n^2\frac{\Lambda}{\lambda}E\tilde{\Delta}(\frac{\partial_n\Phi_1}{|x|})\\
=&\frac{f'}{f}\tilde{\Delta}\Phi_1+\frac{b}{2t^2}\tilde{\Delta}\psi
-(a^{kl}\partial_{kl}a^{ij}_\epsilon+\partial_ka^{kl}\partial_la^{ij}_\epsilon)
(\partial_{ij}\Phi+2\partial_i\Phi\partial_j\Phi)\\
&-(\partial_ka^{kl}a^{ij}_\epsilon+2a^{kl}\partial_ka^{ij}_\epsilon)
(\partial_{ijl}\Phi+4\partial_{il}\Phi\partial_j\Phi)\\
&-a^{kl}a^{ij}_\epsilon(\partial_{ijkl}\Phi+4\partial_{ikl}\Phi\partial_j\Phi+
4\partial_{ik}\Phi\partial_{jl}\Phi)\\
&-8n^2\frac{\Lambda}{\lambda}E[\partial_n^3\Phi_1\frac{a^{nn}}{|x|}+
\partial_{nn}\Phi_1(\frac{\partial_ka^{kn}}{|x|^2}-\frac{2a^{kn}x_k}{|x|^3})\\
&+\partial_n\Phi_1(\frac{3a^{kl}x_k x_l}{|x|^5}-\frac{\partial_ka^{kl}x_l+a^{kk}}{|x|^3})].
\end{split}
\end{eqnarray*}
Next we estimate the terms of $J_6$. By (\ref{indexA}) and (\ref{pc2.2}) we have
\begin{equation*}
\begin{split}
&|a^{ij}_\epsilon|,~|\nabla a^{ij}_\epsilon|~\text{and}~|\nabla^2a^{ij}_\epsilon|~\text{are all bounded};\\
&|\nabla^k\psi|\leq \frac{C}{|x|^{k-2}},~~k=1,2,3,4.
\end{split}
\end{equation*}
Then it is easy to see\\
$$|\tilde{\Delta}\Phi_1|=|\partial_ia^{in}\partial_n\Phi_1+a^{nn}\partial_{nn}\Phi_1|\leq C\partial_n\Phi_1;$$
$$|\tilde{\Delta}\psi|=|\partial_ia^{ij}\partial_j\psi+a^{ij}\partial_{ij}\psi|
\leq C(|\nabla\psi|+|\nabla^2\psi|)\leq C|x|;$$
and
\begin{eqnarray*}
\begin{split}
&(a^{kl}\partial_{kl}a^{ij}_\epsilon+\partial_ka^{kl}\partial_la^{ij}_\epsilon)(\partial_{ij}\Phi+2\partial_i\Phi\partial_j\Phi)\\
\leq& C(|\nabla^2\Phi|+|\nabla\Phi|^2)\\
\leq& C[\partial_{nn}\Phi_1+\frac{|\nabla^2\psi|}{t}+(\partial_n\Phi_1)^2+\frac{|\nabla\psi|^2}{t^2}]\\
\leq& C[(\partial_n\Phi_1)^2+\partial_n\Phi_1+\frac{1}{t}+\frac{|x|^2}{t^2}];
\end{split}
\end{eqnarray*}
also
\begin{eqnarray*}
\begin{split}
&(\partial_ka^{kl}a^{ij}_\epsilon+2a^{kl}\partial_ka^{ij}_\epsilon)(\partial_{ijl}\Phi+
4\partial_{il}\Phi\partial_j\Phi)\\
\leq& C(|\nabla^3\Phi|+|\nabla\Phi|^2+|\nabla^2\Phi|^2)\\
\leq& C[|\partial_n^3\Phi_1|+\frac{|\nabla^3\psi|}{t}+(\partial_n\Phi_1)^2+\frac{|\nabla\psi|^2}{t^2}
+(\partial_{nn}\Phi_1)^2+\frac{|\nabla^2\psi|^2}{t^2}]\\
\leq& C[(\partial_n\Phi_1)^2+\partial_n\Phi_1+\frac{1}{t^2}+\frac{|x|^2}{t^2}];
\end{split}
\end{eqnarray*}

\begin{eqnarray*}
\begin{split}
&a^{kl}a^{ij}_\epsilon(\partial_{ijkl}\Phi+4\partial_{ikl}\Phi\partial_j\Phi+
4\partial_{ik}\Phi\partial_{jl}\Phi)\\
\leq& C(|\nabla^4\Phi|+\sum_{k=1}^3|\nabla^k\Phi|^2)\\
\leq& C[|\partial_n^4\Phi_1|+\frac{|\nabla^4\psi|}{t}+\sum_{k=1}^3(\partial_n^k\Phi_1)^2+
\sum_{k=1}^3\frac{|\nabla^k\psi|^2}{t^2}]\\
\leq& C[(\partial_n\Phi_1)^2+\partial_n\Phi_1+\frac{1}{t^2}+\frac{|x|^2}{t^2}];
\end{split}
\end{eqnarray*}

\begin{eqnarray*}
\begin{split}
&\partial_n^3\Phi_1\frac{a^{nn}}{|x|}+\partial_{nn}\Phi_1(\frac{\partial_ka^{kn}}{|x|^2}-\frac{2a^{kn}x_k}{|x|^3})
+\partial_n\Phi_1(\frac{3a^{kl}x_k x_l}{|x|^5}-\frac{\partial_ka^{kl}x_l+a^{kk}}{|x|^3})\\
\leq&C(|\partial_n^3\Phi_1|+\partial_{nn}\Phi_1+\partial_n\Phi_1)\\
\leq&C\partial_n\Phi_1.
\end{split}
\end{eqnarray*}
Combining them together, we have
\begin{equation}\label{pc2.14}
J_6\geq C\frac{f'}{f}\partial_n\Phi_1-C[(\partial_n\Phi_1)^2+\partial_n\Phi_1+\frac{1}{t^2}
+\frac{|x|^2}{t^2}].
\end{equation}
Combining (\ref{pc2.6}), (\ref{pc2.7}), (\ref{pc2.8}),
(\ref{pc2.12}), (\ref{pc2.13}) and (\ref{pc2.14}), we have
\begin{eqnarray*}
\begin{split}
M_2\geq&2[(\alpha-1)\lambda^2-(16n^2\frac{\Lambda}{\lambda}+8n^2+4n)\Lambda E]\frac{(\partial_n\Phi_1)^3}{|x|}\\
&+[-(4\lambda-(16n^2\frac{\Lambda}{\lambda}+2n)E)\frac{f'}{f}-\frac{2d\Lambda+C}{t}](\partial_n\Phi_1)^2\\
&+[\partial_{tt}\Phi_1+\frac{2d}{t}\partial_t\Phi_1+C\frac{f'}{f}\partial_n\Phi_1-
C\partial_n\Phi_1+C(\frac{f'}{tf}-\frac{1}{t^2})\Phi_1]\\
&+[\frac{b\psi}{t^3}(\frac{d}{2}-C)-\frac{C}{t^3}|x|^2-\frac{db^2\Lambda+C}{2t^3}|\nabla\psi|^2]+
\frac{Kd}{4t^3}-\frac{2d^2+C}{t^2}.
\end{split}
\end{eqnarray*}
Next we estimate the terms of the right side of the above inequality.
We always choose $d$ large enough.\\

Notice that
$(16n^2\frac{\Lambda}{\lambda}+2n)E<16n^2(\frac{\Lambda}{\lambda}+1)E_0<\lambda$,
then
\begin{eqnarray*}
\begin{split}
&-(4\lambda-(16n^2\frac{\Lambda}{\lambda}+2n)E)\frac{f'}{f}-\frac{2d\Lambda+C}{t}\\
\geq&-3\lambda\frac{f'}{f}-\frac{3d\Lambda}{t}\\
=&\frac{3\lambda K}{t(1-t^K)}-\frac{3d\Lambda}{t}\\
\geq&\frac{3\lambda K-3d\Lambda}{t}\\
\geq&0.
\end{split}
\end{eqnarray*}
And
\begin{eqnarray*}
\begin{split}
&\partial_{tt}\Phi_1+\frac{2d}{t}\partial_t\Phi_1+C\frac{f'}{f}\partial_n\Phi_1-C\partial_n\Phi_1+
C(\frac{f'}{tf}-\frac{1}{t^2})\Phi_1\\
\geq&\Phi_1[\frac{f''}{f}+\frac{2d}{t}\frac{f'}{f}+C\frac{f'}{f}-C+C(\frac{f'}{tf}-\frac{1}{t^2})]\\
\geq&\Phi_1[\frac{f''}{f}+\frac{3d}{t}\frac{f'}{f}-\frac{d}{t^2}]\\
=&\Phi_1[\frac{K(K+1-3d)}{t^2(1-t^K)}-\frac{d}{t^2}]\\
\geq&\Phi_1\frac{K(K+1-3d)-d}{t^2}\\
\geq&0.
\end{split}
\end{eqnarray*}
By (\ref{pc2.2})
$$\psi\geq\frac{|x|^2}{2},~~~|\nabla \psi|\leq 4(\frac{\Lambda}{\lambda}+1)^2|x|,$$
then
\begin{eqnarray*}
\begin{split}
&\frac{b\psi}{t^3}(\frac{d}{2}-C)-\frac{C}{t^3}|x|^2-\frac{db^2\Lambda+C}{2t^3}|\nabla\psi|^2\\
\geq&(\frac{db}{4}-C)\frac{|x|^2}{t^3}-C\frac{|x|^2}{t^3}-8(db^2\Lambda+C)
(\frac{\Lambda}{\lambda}+1)^4\frac{|x|^2}{t^3}\\
=&\frac{|x|^2}{t^3}[\frac{db}{4}-8db^2\Lambda(\frac{\Lambda}{\lambda}+1)^4-C]\\
=&\frac{|x|^2}{t^3}[db(\frac{1}{4}-8b\Lambda(\frac{\Lambda}{\lambda}+1)^4)-C].
\end{split}
\end{eqnarray*}
Since
$$b=\frac{1}{64\Lambda(\frac{\Lambda}{\lambda}+1)^4},$$
then
$$
\frac{b\psi}{t^3}(\frac{d}{2}-C)-\frac{C}{t^3}|x|^2-\frac{db^2\Lambda+C}{2t^3}|\nabla\psi|^2\geq
\frac{|x|^2}{t^3}(\frac{db}{8}-C)
\geq\frac{db|x|^2}{16t^3}.
$$
Finally,
$$\frac{Kd}{4t^3}-\frac{2d^2+C}{t^2}\geq\frac{Kd}{4t^3}-\frac{3d^2}{t^2}\geq\frac{Kd-12d^2}{4t^3}\geq
\frac{1}{t^3}.$$ Combining them together, we have
$$M_2\geq2[(\alpha-1)\lambda^2-(16n^2\frac{\Lambda}{\lambda}+8n^2+4n)\Lambda E]\frac{(\partial_n\Phi_1)^3}{|x|}
+\frac{bd|x|^2}{16t^3}+\frac{1}{t^3}.$$

\emph{Estimate of $|\nabla(F-F_0)|$}.\\

Recall that
$$F-F_0=2(a^{ij}_\epsilon-a^{ij})(\partial_{ij}\Phi+2\partial_i\Phi\partial_j\Phi),$$
then
\begin{eqnarray*}
\begin{split}
|\nabla(F-F_0)|=&2|(\nabla a^{ij}_\epsilon-\nabla a^{ij})(\partial_{ij}\Phi+2\partial_i\Phi\partial_j\Phi)\\
&+(a^{ij}_\epsilon-a^{ij})(\nabla\partial_{ij}\Phi+4\partial_i\Phi\nabla\partial_j\Phi)|\\
\leq&2|\nabla a^{ij}_\epsilon-\nabla a^{ij}|(|\partial_{ij}\Phi|+|\partial_i\Phi|^2+|\partial_j\Phi|^2)\\
&+2|a^{ij}_\epsilon-a^{ij}|(|\nabla\partial_{ij}\Phi|+2|\partial_i\Phi|^2+2|\nabla\partial_j\Phi|^2).
\end{split}
\end{eqnarray*}
Since $|x|\geq1$ in $Q$, by (\ref{indexA}) we have
$$|\nabla a^{ij}_\epsilon|\leq\frac{2E}{|x|},~~|a^{ij}_\epsilon-a^{ij}|\leq\frac{E}{|x|},$$
then
\begin{eqnarray*}
\begin{split}
|\nabla(F-F_0)|\leq&\frac{6E}{|x|}(|\partial_{nn}\Phi_1|+\frac{C}{t}|\nabla^2\psi|+2n|\nabla\Phi|^2)\\
&+\frac{2E}{|x|}(|\partial_n^3\Phi_1|+\frac{C}{t}|\nabla^3\psi|+2n|\nabla\Phi|^2+2n|\nabla^2\Phi|^2)\\
\leq&\frac{E}{|x|}[8\partial_n\Phi_1+\frac{C}{t}(|\nabla^2\psi|+|\nabla^3\psi|)+16|\nabla\Phi|^2+
4n|\nabla^2\Phi|^2].
\end{split}
\end{eqnarray*}
Next we estimate the terms of the right side of the above inequality.\\

By the inequality $(a+b)^2\leq\frac{3}{2}a^2+3b^2$, we have
$$|\nabla\Phi|^2=|\nabla\Phi_1+\nabla\Phi_2|^2\leq\frac{3}{2}|\nabla\Phi_1|^2+
3|\nabla\Phi_2|^2\leq\frac{3}{2}(\partial_n\Phi_1)^2+\frac{C|x|^2}{t^2};$$
$$|\nabla^2\Phi|^2=|\nabla^2\Phi_1+\nabla^2\Phi_2|^2\leq
\frac{3}{2}|\nabla^2\Phi_1|^2+3|\nabla^2\Phi_2|^2\leq\frac{3}{2}(\partial_{nn}\Phi_1)^2+\frac{C}{t^2}.$$
Then
\begin{eqnarray*}
\begin{split}
|\nabla(F-F_0)|\leq&\frac{E}{|x|}[8\partial_n\Phi_1+24n(\partial_n\Phi_1)^2+
6n(\partial_{nn}\Phi_1)^2+\frac{C|x|^2}{t^2}]\\
\leq&\frac{E}{|x|}[8\partial_n\Phi_1+30n(\partial_n\Phi_1)^2+\frac{C|x|^2}{t^2}].
\end{split}
\end{eqnarray*}
Because
$$8\partial_n\Phi_1\leq2n(\partial_n\Phi_1)^2+\frac{8}{n}\leq2n(\partial_n\Phi_1)^2+\frac{8|x|^2}{t^2},$$
then
$$
|\nabla(F-F_0)|\leq\frac{E}{|x|}[32n(\partial_n\Phi_1)^2+\frac{C|x|^2}{t^2}]
=32nE[\frac{(\partial_n\Phi_1)^2}{|x|}+\frac{C|x|}{t^2}].
$$
Thus we proved Lemma \ref{estimates2}.

\section{Appendix}

\emph{Appendix A: The properties of $\{a^{ij}_\epsilon\}$.}\\

$a^{ij}_\epsilon(x,t)=\int_{\mathbb{R}^n}a^{ij}(x-y,t)\phi_\epsilon(y)dy$,
where $\phi$ is a mollifier and $\epsilon=\frac{1}{2}$, and
$$
(x,t)\in \left\{
           \begin{array}{ll}
             \mathbb{R}^n\times(0,2), & \hbox{under the assumptions of Proposition \ref{Prop-C1};} \\
             (\mathbb{R}^n_++e_n)\times(0,1), & \hbox{under the assumptions of Proposition \ref{Prop-C2}.}
           \end{array}
         \right.
$$
Then $\{a^{ij}_\epsilon\}$ satisfy:
\begin{eqnarray*}
\begin{split}
&i)~\lambda|\xi|^2\leq a^{ij}_\epsilon(x,t)\xi_i\xi_j\leq \Lambda|\xi|^2, ~\forall \xi\in\mathbb{R}^n;\\
&ii)~|\nabla a^{ij}_\epsilon(x,t)|\leq M;~~|\nabla a^{ij}_\epsilon(x,t)|\leq\frac{2E}{|x|}~when~|x|\geq1;\\
&iii)~|a^{ij}_\epsilon(x,t)-a^{ij}(x,t)|\leq 2\Lambda;
~~|a^{ij}_\epsilon(x,t)-a^{ij}(x,t)|\leq\frac{E}{|x|}~when~|x|\geq1;\\
&iv)~|\partial_{kl}a^{ij}_\epsilon(x,t)|\leq c(n)M;~~|\partial_{kl}a^{ij}_\epsilon(x,t)|\leq \frac{c(n)E}{|x|}~when~|x|\geq1.
\end{split}
\end{eqnarray*}
\emph{Proof}.\\
$i)$ Obvious.\\
$ii)$
$$|\nabla a^{ij}_\epsilon(x,t)|\leq\int_{\mathbb{R}^n}|\nabla a^{ij}(x-y,t)|\phi_\epsilon(y)dy
\leq M\int_{\mathbb{R}^n}\phi_\epsilon(y)dy=M,$$
and when $|x|\geq1$,
$$
|\nabla a^{ij}_\epsilon(x,t)|\leq\int_{\mathbb{R}^n}|\nabla a^{ij}(x-y,t)|\phi_\epsilon(y)dy
\leq \int_{\mathbb{R}^n}\frac{E}{|x-y|}\phi_\epsilon(y)dy
\leq \int_{\mathbb{R}^n}\frac{E}{|x|-\frac{1}{2}}\phi_\epsilon(y)dy
\leq\frac{2E}{|x|}.
$$
$iii)$ The first part is obvious. We only need to prove the second
one.
\begin{eqnarray*}
\begin{split}
|a^{ij}_\epsilon(x,t)-a^{ij}(x,t)|&\leq\int_{\mathbb{R}^n}|a^{ij}(x-y,t)-a^{ij}(x,t)|\phi_\epsilon(y)dy\\
&\leq\int_{\mathbb{R}^n}|\nabla a^{ij}(x-\theta y,t)||y|\phi_\epsilon(y)dy, ~~~~~~~(0<\theta<1)\\
\end{split}
\end{eqnarray*}
and when $|x|\geq1$,
$$
|a^{ij}_\epsilon(x,t)-a^{ij}(x,t)|\leq\int_{\mathbb{R}^n}\frac{E}{2|x-\theta y|}\phi_\epsilon(y)dy
\leq\int_{\mathbb{R}^n}\frac{E}{2(|x|-\frac{1}{2})}\phi_\epsilon(y)dy
\leq\frac{E}{|x|}.
$$
$iv)$
\begin{eqnarray*}
\begin{split}
|\partial_{kl}a^{ij}_\epsilon(x,t)|&\leq\int_{\mathbb{R}^n}|
\partial_ka^{ij}(x-y,t)||\partial_l\phi_\epsilon(y)|dy\\
&\leq\epsilon^{-n-1}\int_{\mathbb{R}^n}|\partial_ka^{ij}(x-y,t)||(\partial_l\phi)(\frac{y}{\epsilon})|dy\\
&\leq\frac{M}{\epsilon}||\partial_l\phi||_1\\
&\leq 2||\nabla\phi||_1M,
\end{split}
\end{eqnarray*}
and when $|x|\geq1$,
\begin{eqnarray*}
\begin{split}
|\partial_{kl}a^{ij}_\epsilon(x,t)|&\leq\epsilon^{-n-1}\int_{\mathbb{R}^n}|
\partial_ka^{ij}(x-y,t)||(\partial_l\phi)(\frac{y}{\epsilon})|dy\\
&\leq\epsilon^{-n-1}\int_{\mathbb{R}^n}\frac{E}{|x-y|}|(\partial_l\phi)(\frac{y}{\epsilon})|dy\\
&\leq \frac{2E}{\epsilon|x|}||\partial_l\phi||_1\\
&\leq \frac{4||\nabla\phi||_1E}{|x|}.
\end{split}
\end{eqnarray*}

\emph{Appendix B: The properties of $\psi$.}

Recall that
$\psi(x)=|x|^2-2\frac{\Lambda}{\lambda}|x|x_n+2(\frac{\Lambda}{\lambda})^2x_n^2$.
Then $\psi$ satisfies
\begin{eqnarray*}
\begin{split}
&i)~\psi\geq\frac{|x|^2}{2};\\
&ii)~D^2\psi\leq C(\frac{\Lambda}{\lambda})I_n;\\
&iii)~|\nabla \psi|\leq 4(\frac{\Lambda}{\lambda}+1)^2|x|, ~|\nabla^k\psi|\leq
\frac{C(n,\frac{\Lambda}{\lambda})}{|x|^{k-2}},k=2,3,4;\\
&iv)~a^{ni}\partial_i\psi\leq C(\Lambda,\lambda)x_n.
\end{split}
\end{eqnarray*}
\emph{Proof}.
Because $i)$, the first part of $iii)$ and $iv)$ play more important role in the proof of
Proposition \ref{Prop-C2}, we just prove them. The others can be proved by direct calculations.\\
In fact
$$\psi=\frac{|x|^2}{2}+\frac{1}{2}(|x|-2\frac{\Lambda}{\lambda}x_n)^2\geq\frac{|x|^2}{2},$$
and
$$
\partial_i\psi=\left\{
                    \begin{array}{ll}
                      2(1-\frac{\Lambda}{\lambda}\frac{x_n}{|x|})x_i, & \hbox{if $1\leq i\leq n-1$;} \\
                      -2\frac{\Lambda}{\lambda}|x|-2\frac{\Lambda}{\lambda}\frac{x_n}{|x|}x_n+
                      [4(\frac{\Lambda}{\lambda})^2+2]x_n, & \hbox{if $i=n$.}
                    \end{array}
                  \right.
$$
Then
\begin{eqnarray*}
\begin{split}
|\nabla\psi|&\leq2|1-\frac{\Lambda}{\lambda}\frac{x_n}{|x|}||x'|
+|-2\frac{\Lambda}{\lambda}|x|-2\frac{\Lambda}{\lambda}\frac{x_n}{|x|}x_n+[4(\frac{\Lambda}{\lambda})^2+2]x_n|\\
&\leq2(1+\frac{\Lambda}{\lambda})|x|+[4(\frac{\Lambda}{\lambda})^2+4\frac{\Lambda}{\lambda}+2]|x|\\
&\leq4(\frac{\Lambda}{\lambda}+1)^2|x|.
\end{split}
\end{eqnarray*}
Direct calculations show that
$$a^{ni}\partial_i\psi=2\sum_{i=1}^{n-1}a^{ni}x_i-2\frac{\Lambda}{\lambda}a^{nn}|x|-
2\frac{\Lambda}{\lambda}\frac{x_n}{|x|}\sum_{i=1}^na^{ni}x_i+[4(\frac{\Lambda}{\lambda})^2+2]a^{nn}x_n.$$
By the Cauchy inequality we have
$$|\sum_{i=1}^na^{ni}x_i|\leq\sqrt{\sum_{i=1}^n(a^{ni})^2\sum_{i=1}^nx_i^2}\leq\Lambda|x|.$$
Similarly, we have
$$|\sum_{i=1}^{n-1}a^{ni}x_i|\leq\Lambda|x'|.$$
Since $$\lambda\leq a^{nn}\leq\Lambda,$$
then
\begin{eqnarray*}
\begin{split}
a^{ni}\partial_i\psi\leq&2\Lambda|x'|-2\Lambda|x|+2\frac{\Lambda^2}{\lambda}x_n+
[4(\frac{\Lambda}{\lambda})^2+2]\Lambda x_n\\
\leq&[4(\frac{\Lambda}{\lambda})^2+2\frac{\Lambda}{\lambda}+2]\Lambda x_n.
\end{split}
\end{eqnarray*}

\bigskip

%
%\noindent {\bf Acknowledgments.}
%Authors thank helpful discussions with Professors.
%

\vspace{1cm}

{\small}

\noindent

{\bf Jie Wu}\\
Institute of Mathematics, Academy of Mathematics and Systems Science, Chinese Academy of Sciences, Beijing 100190, P. R. China.\\
Email address: {\bf jackwu@amss.ac.cn}\\

{\bf Liqun Zhang}\\
Hua Loo-Keng Key Laboratory of Mathematics, Chinese Academy of Sciences, Beijing 100190, P. R. China. \\
Email address: {\bf lqzhang@amss.ac.cn}

\end{document}